\documentclass[12pt]{article}
\usepackage{amsfonts}
\usepackage{amssymb}
\usepackage{polski}
\input{tcilatex}
\begin{document}

\begin{center}
\textbf{Category-measure duality: convexity, mid-point convexity and Berz
sublinearity\\[0pt]
by \\[0pt]
N. H. Bingham and A. J. Ostaszewski}\\[0pt]

\bigskip
\end{center}

\noindent \textbf{Abstract.}

Category-measure duality concerns applications of Baire-category methods
that have measure-theoretic analogues. The set-theoretic axiom needed in
connection with the Baire category theorem is the Axiom of Dependent Choice
DC rather than the Axiom of Choice AC. Berz used the Hahn-Banach Theorem
over $\mathbb{Q}$ to prove that the graph of a measurable sublinear function
that is $\mathbb{Q}_{+}$-homogeneous consists of two half-lines through the
origin. We give a category form of the Berz theorem. Our proof is simpler
than that of the classical measure-theoretic Berz theorem, our result
contains Berz's theorem rather than simply being an analogue of it, and we
use only DC rather than AC. Furthermore, the category form easily
generalizes: the graph of a Baire sublinear function defined on a Banach
space is a cone. The results are seen to be of automatic-continuity type. We
use Christensen Haar null sets to extend the category approach beyond the
locally compact setting where Haar measure exists. We extend Berz's result
from Euclidean to Banach spaces, and beyond. Passing from sublinearity to
convexity, we extend the Bernstein-Doetsch theorem and related continuity
results, allowing our conditions to be `local' -- holding off some
exceptional set.

\bigskip



\section{Introduction}

The Berz theorem of our title is his characterization of a function $S:%
\mathbb{R}\rightarrow \mathbb{R}$ which is \textit{sublinear,} that is -- it
is \textit{subadditive }([HilP, Ch. 3], [Ros]):%
\[
S(u+v)\leq S(u)+S(v),
\]%
and \textit{homogeneous} with respect to non-negative integer scaling.
Following Berz [Ber], we call $S$ \textit{sublinear} \textit{on a set} $%
\Sigma $ if $S$ is subadditive and
\[
S(nx)=nS(x)\text{ for }x\in \Sigma ,n=0,1,2,...,
\]%
equivalently, if $\Sigma $ is closed under non-negative rational scaling,%
\[
S(qx)=qS(x)\text{ for }x\in \Sigma ,q\in \mathbb{Q}_{+}=\mathbb{Q\cap
\lbrack }0,\infty );
\]%
in words, $S$ is positively $\mathbb{Q}$-homogeneous on $\Sigma $ and $S(0)=0
$. An important class of functions with these two properties but with a more
general domain occurs in mathematical finance -- the \textit{coherent risk
measures} introduced by Artzner et al. [ArtDEH] (cf. \S 6.5); for textbook
treatments see [McNFE], [FolS, 4.1]. In \S 4 we characterize such functions
in the category setting when the domain is a Banach space. Working in a
locally convex Fr\'{e}chet space and under various axiomatic assumptions
Ajtai [Ajt], Wright [Wri], and Garnir [Gar], motivated by semi-norm
considerations, study the continuity of a subadditive function $S$ with the
property $S(2x)=2S(x)$.

Berz used the Hahn-Banach Theorem over $\mathbb{Q}$ to prove that the graph
of a (Lebesgue) measurable sublinear function consists of two half-lines
through the origin ([Kuc, \S 16.4,5]; cf. [BinO1]). Recall that in a
topological space $X,$ a subset $H$ is \textit{Baire} (has the \textit{Baire
property, BP}) if $H=(V\backslash M_{V})\cup M_{H}$ for some open set $V$
and meagre sets $M_{V},M_{H}$ in the sense of the topology on $X$; similarly
a function $f:X\rightarrow \mathbb{R}$ is \textit{Baire} if preimages of
(Euclidean) open subsets of $\mathbb{R}$ are Baire subsets in the topology
of $X$. Our first result is the Baire version of Berz's theorem on the line.
Below $\mathbb{R}_{\pm }$ denotes the non-negative and non-positive
half-lines.

\bigskip

\noindent \textbf{Theorem 1B }(ZF+DC)\textbf{.} \textit{For} $S:\mathbb{R}%
\rightarrow \mathbb{R}$ \textit{sublinear and Baire, there are} $c^{\pm }\in
\mathbb{R}$ \textit{such that}%
\[
S(x)=c^{\pm }x,\text{ for }x\in \mathbb{R}_{\pm }.
\]

As we shall see in \S 3, Theorem 1B implies the classical Berz theorem as a
corollary:

\bigskip

\noindent \textbf{Theorem 1M }(ZF+DC, containing Berz [Ber] with AC)\textbf{.%
} \textit{For} $S:\mathbb{R}\rightarrow \mathbb{R}$ \textit{sublinear and
measurable, there are} $c^{\pm }\in \mathbb{R}$ \textit{such that}%
\[
S(x)=c^{\pm }x,\text{ for }x\in \mathbb{R}_{\pm }.
\]

Theorems 1B and 1M may be combined, into `Theorem 1(B+M)', say. Following
necessary topological preliminaries (Lemma S, Theorem BL; Steinhaus-Weil
property) in \S 2, the two cases are proved together in \S 3
bi-topologically, by switching between the two relevant density topologies
of \S 2 ([BinO6,10,15], [Ost2]). Here we also prove Theorem 2 (local
boundedness for subadditive functions) and Theorem 3, the corresponding
continuity result. We introduce universal measurability (used in \S 4 in
defining Christensen's notion of Haar null sets -- in contexts where there
may be no Haar measure -- [Chr1,2]), and use this to note a variant on
Theorem 2, Theorem 2H (`H for Haar').

The sector between the lines $c^{\pm }x$ in the upper half-plane is a
two-dimensional cone. This suggests the generalization to Banach spaces that
we prove in \S 4 (Theorems 4B, 4M, 4F -- `F for F-space').

The results above for the Baire/measurable functions on $\mathbb{R}$ are to
be expected: they follow from the classical Bernstein-Doetsch continuity
theorem for locally bounded mid-point convex functions on normed vector
spaces, to which we turn in \S 5 (see e.g. [Kuc, 6.4.2] quoted for $\mathbb{R%
}^{d}$, but its third proof there applies more generally, as does Theorem B\
below, also originally for $\mathbb{R}^{d}$; see also [HarLP, Ch. III]),
once one proves their local boundedness (\S 3, Th. 2), since a sublinear
function is necessarily mid-point convex. Indeed, by $\mathbb{Q}$%
-homogeneity and subadditivity,%
\[
f\left( \frac{1}{2}(x+y)\right) =\frac{1}{2}\left( f(x+y)\right) \leq \frac{1%
}{2}\left( f(x)+f(y)\right) .
\]%
We handle the Berz sublinear case first (in \S 3), as the arguments are
simpler, and turn to mid-convexity matters in \S 5, where we prove the
following two results (for topological and convexity terminology see
respectively \S 2 and 5).

\bigskip

\noindent \textbf{Theorem M }(Mehdi's Theorem, [Meh, Th. 4]; cf. [Wri]).
\textit{For a Banach space }$X,$ \textit{if }$S:X\rightarrow \mathbb{R}$%
\textit{\ is mid-point convex and Baire,} \textit{then }$S\ $\textit{is
continuous.}

\bigskip

\noindent \textbf{Theorem FS }(cf. [FisS])\textbf{. }\textit{For a Banach
space }$X,$ \textit{if }$S:X\rightarrow \mathbb{R}$\textit{\ is mid-point
convex and universally measurable,} \textit{then }$S\ $\textit{is continuous.%
}

\bigskip

For the Banach context both there and in \S 4, we rely on the following
dichotomy result, Theorem B, especially on its second assertion, which
together with an associated Corollary B in \S 4 (on boundedness), enables
passage from a general Banach space to a separable one (wherein the
Christensen theory of Haar null sets is available). See [Blu] and Appendix 2
of the arXiv version of this paper.

\bigskip

\noindent \textbf{Theorem B (Blumberg's Dichotomy Theorem}, [Blu, Th. 1];
cf. [Sie2]\textbf{).} \textit{For }$X$\textit{\ any normed vector space and }%
$S:X\rightarrow \mathbb{R}$ \textit{mid-point convex: either }$S$\textit{\
is not continuous at }$x_{0}\in X,$\textit{\ or }$S(x_{n})$\textit{\ is
unbounded above for some sequence }$x_{n}$\textit{\ with limit }$x_{0}$%
\textit{.}

\textit{In particular, for }$X$\textit{\ a Banach space, if for any closed
separable subspace }$B\subseteq X$\textit{\ the restriction }$S|B$\textit{\
is continuous (for instance }$S|B$ \textit{is locally bounded on }$B$\textit{%
), then }$S$\textit{\ is continuous.}

\bigskip

In \S 5 we switch to a form of mid-convexity that is assumed to hold only on
a co-meagre or co-null set (so on an open set of a \textit{density topology}
-- see \S 2); we term this \textit{weak mid-point convexity,} and show in
particular that a Baire/measurable weakly mid-point convex function is
continuous and convex. It follows that the Berz theorems are true under the
hypothesis of \textit{weak sublinearity} (sublinearity on a co-meagre or
co-null set); however, we leave open the possibility of a direct proof along
the lines of \S 3 and also the question whether a Bernstein-Doetsch
dichotomy holds -- that a weakly mid-point convex/sublinear function is
either everywhere continuous or nowhere continuous. We close in \S 6 with
some complements.

Theorem 1B (under the usual tacit assumption ZF+AC) was given in [BinO13,
Th. 5]. The results imply the classical results that Baire/measurable
additive functions are linear (see [BinO9] for historical background);
indeed, an additive function $A(.)$ is sublinear and $A(-x)=-A(x),$ so $%
c^{+}=-c^{-}.$

The primacy of category within category-measure duality is one of our two
main themes here. This is something we have emphasised before [BinO6,9,10];
Oxtoby [Oxt] calls this measure-category duality, but from a different
viewpoint -- he has no need of Steinhaus's theorem (cf. [Ost2]), which is
crucial for us. Our second main theme, new here, is AC\ versus DC. As so
much of the extensive relevant background is still somewhat scattered, we
summarize what we need in detail in Appendix 1 (which has its own separate
references). This may be omitted by the expert (or uninterested) reader, and
so is included only in the fuller arXiv version of this paper.

Without further comment, we work with ZF+DC, rather than ZF+AC, throughout
the paper. It is natural that DC should dominate here. For, DC suffices for
the \textit{common parts} of the Baire category and Lebesgue measure cases:
for the first, see Blair [Bla], and for the second, see Solovay (Appendix
1.3; [Solo2, p. 25]). For the \textit{contrasts} -- or `wedges' -- between
them, see Appendix 1.5. It is here that further set-theoretic assumptions
become crucial; in brief, \textit{measure theory needs stronger assumptions}.

\section{Topological preliminaries: Steinhaus-Weil property}

Fundamental for our purposes is the \textit{Steinhaus-Weil property}%
\footnote{%
Initially, as in the Steinhaus-Piccard-Pettis context, this concerns $%
\mathbb{R}$; the wider context is due to Weil and concerns (Haar)
measurability in locally compact groups [Wei, p. 50], cf. [GroE2]. These
distinctions blur in our bitopological context.} [BinO14,15] -- that the
difference set $A-A$ has non-empty interior for $A$ any non-negligible set
with the Baire property, briefly: \textit{Baire set} -- as opposed to Baire
\textit{topology}. We focus on\textit{\ Baire }topological spaces on which
the Steinhaus-Weil theorem holds. (See [Sole2, Remark to Th. 6.1] for
failure of the Steinhaus-Weil property in a group; cf. [Kom] and [RosS] for
extensions of this property.) This is just what is needed to make the
infinite combinatorics used in our proofs work.

Call $a$ an (outer) \textit{Lebesgue-density} point of a set $A$ if $%
\lim_{\delta \downarrow 0}|A\cap (a-\delta ,a+\delta )|/2\delta =1,$ where $%
|S|$ is the outer measure of $S;$ the \textit{Lebesgue density theorem}
asserts that almost all points of a set are density points. (On this point
the survey [Bru] is a classic. For further background see [Wil] and
literature cited there; cf. the recent [BinO14].) By analogy, say that $a$
is a \textit{Baire-density} point of $A$ if $V\backslash A$ is meagre, for
some open neighbourhood $V$ of $a;$ if $A$ is Baire, then it is immediate
from the BP that, except for a meagre set, all points of $A$ are Baire
density points. Each of the category and measure notions of density defines
a \textit{density topology} (denoted respectively $\mathcal{D}_{\mathcal{B}}/%
\mathcal{D}_{\mathcal{L}}$ -- with $\mathcal{L}$ denoting Lebesgue
measurable sets), in which a set $W$ is \textit{density-open} if all its
points are category/measure density points of $W,$ the latter case
introduced by Goffman and his collaborators -- see [GofNN] and [GofW]. Both
refine the usual Euclidean topology, $\mathcal{E}$; see [BinO14] for
properties common to both topologies. We call meagre/null sets \textit{%
negligible}, and say that \textit{quasi all} points of a set have a property
if, but for a negligible subset, all have the property. These negligible
sets form a $\sigma $-ideal; see Fremlin [Fre2], Lukes et al. [LukMZ], Wilczy%
\'{n}ski [Wil], [BinO14] for background, and also [BinO5] and [Ost1]. Below
(for use in \S 4) we consider a further $\sigma $-ideal: the \textit{left
Haar null sets} (equivalently: \textit{Haar null}) of a Banach space and by
extension use the same language of negligibles there. The corresponding
density topologies may also be studied via the \textit{Hashimoto} topologies
(cf. [Has], [BalR], [LukMZ, 1C]), obtained by declaring as basic open the
sets of the form $U\backslash N$ with $U\in \mathcal{E}$ and $N$ the
appropriate negligible. (That these sets, even under DC, form a topology
follows from $\mathcal{E}$ being second countable -- cf. [JanH, 4.2] and
[BinO14].)

The definition above of a Baire-density point may of course be repeated
verbatim in the context of any topology $\mathcal{T}$ on any set $X$ by
referring to $\mathcal{B}(\mathcal{T}),$ the Baire sets of $\mathcal{T}$ .
In particular, working with $\mathcal{T}=\mathcal{D}_{\mathcal{L}}$ in place
of $\mathcal{E}$ we obtain a topology $\mathcal{D}_{\mathcal{B}}(\mathcal{D}%
_{\mathcal{L}}).$ Since $\mathcal{B}(\mathcal{D}_{\mathcal{L}})=\mathcal{L}$
(see [Kec, 17.47] and [BinO6]),%
\[
\mathcal{D}_{\mathcal{B}}=\mathcal{D}_{\mathcal{B}}(\mathcal{E}),\qquad
\mathcal{D}_{\mathcal{L}}=\mathcal{D}_{\mathcal{B}}(\mathcal{D}_{\mathcal{L}%
}).
\]

\bigskip

\noindent \textbf{Lemma S (Multiplicative Sierpi\'{n}ski Lemma}; [BinO5,
Lemma S], cf. [Sie1])\textbf{.} \textit{For }$A,B$\textit{\ Baire/measurable
in }$(0,\infty )$ \textit{with respective density points (in the
category/measure sense) }$a,b,$\textit{\ then for }$n=1,2,...$\textit{\
there exist positive rationals }$q_{n}$ \textit{and points }$a_{n},b_{n}$
\textit{converging (metrically) to }$a,b$ \textit{through }$A,B$ \textit{%
respectively such that }$b_{n}=q_{n}a_{n}.$

\bigskip

\noindent \textbf{Proof.} For $n=1,2,...$ and the consecutive values $%
\varepsilon =1/n,$ the sets $B_{\varepsilon }(a)\cap A$ and $B_{\varepsilon
}(b)\cap B$ are Baire/measurable non-negligible. So by Steinhaus's theorem
(see e.g.[Kuc, \S 3.7], [BinGT, Th. 1.1.1]; cf. [BinO9]), the set $[B\cap
B_{\varepsilon }(b)]\cdot \lbrack A\cap B_{\varepsilon }(a)]^{-1}$ contains
interior points, and so in particular a rational point $q_{n}.$ Thus for
some $a_{n}\in B_{\varepsilon }(a)\cap A$ and $b_{n}\in B_{\varepsilon
}(b)\cap B$ we have $q_{n}=b_{n}a_{n}^{-1}>0,$ and as $|a-a_{n}|<1/n$ and $%
|b-b_{n}|<1/n,$ $a_{n}\rightarrow a,b_{n}\rightarrow b$. $\square $

\bigskip

\noindent \textbf{Remark. }The result above is a consequence of the
Steinhaus-Weil property regarded as a corollary of the Category Interior
Theorem ([BinO7, Th. 4.4]; cf. [GroE1,2]). The latter, applied to the
topology $\mathcal{D}$ that is either of the above two topologies $\mathcal{D%
}_{\mathcal{B}}/\mathcal{D}_{\mathcal{L}},$ asserts that $U-V\ $or $UV^{-1}$
is an $\mathcal{E}$-open nhd (of the relevant neutral element) for $U,V$
open under $\mathcal{D}$, since $\mathcal{D}$ is a shift-invariant Baire
topology satisfying the Weak Category Convergence condition of [BinO6] for
either of the shift actions $x\mapsto x+a,$ $x\mapsto xa$. The Category
Interior Theorem in turn follows from the Category Embedding Theorem
([BinO6]; cf. [MilO]). Now $a\in A^{o}:=\mathrm{int}_{\mathcal{D}}(A),$ $%
b\in B^{o},$ as $a$ and $b$ are respectively density points of $A$ and $B,$
and $B_{\varepsilon }(a)\cap A^{o}$ and $B_{\varepsilon }(b)\cap B^{o}$ are
in $\mathcal{D}$, as $\mathcal{D}$ refines $\mathcal{E}$.

\bigskip

\noindent \textbf{Definition. }Say that $f:X\rightarrow \mathbb{R}$ is
\textit{quasi }$\sigma $\textit{-continuous} if $X$ contains a $\mathcal{B}$%
-set $\Sigma ^{+}$\ which is quasi all of $X$\ and an increasing
decomposition $\Sigma ^{+}:=\bigcup_{m=0}^{\infty }\Sigma _{m}$ into $%
\mathcal{B}$-sets $\Sigma _{m}$\ such that each $f|\Sigma _{m}$\ is
continuous.

\bigskip

Separability is a natural condition in the next result -- see the closing
comments in [Zak].

\bigskip

\noindent \textbf{Theorem BL (Baire Continuity Theorem }[BinO8, Th. 11.8]%
\textbf{; Baire-Luzin Theorem; }cf. [Hal], end of Section 55, [Zak, Th.
II]). \textit{For a separable Banach space, if }$f:X\rightarrow \mathbb{R}$%
\textit{\ is Baire, or measurable with respect to a regular }$\sigma $%
\textit{-finite measure, then }$f$\textit{\ is quasi }$\sigma $\textit{%
-continuous, with the sets }$\Sigma _{m}$\textit{\ in the Baire case being
in }$\mathcal{D}_{\mathcal{B}}$\textit{. Furthermore, for }$X=\mathbb{R}$
\textit{under Lebesgue measure the sets }$\Sigma _{m}$ \textit{may likewise
be taken in }$\mathcal{D}_{\mathcal{L}}=\mathcal{D}_{\mathcal{B}}(\mathcal{D}%
_{\mathcal{L}}).$

\bigskip

\noindent \textbf{Remarks. }1. In the category case, with $\Sigma
_{m}=\Sigma _{0}$ for all $m$ and $\Sigma _{0}$ co-meagre, this is Baire's
Theorem ([Oxt, Th. 8.1]). In the Lebesgue measure case this is a useful form
of Luzin's Theorem formulated in [BinO4]. The extension to a regular (i.e. $%
\mathcal{G}$-outer regular) $\sigma $-finite measure may be made via
Egoroff's Theorem (cf. [Hal, \S 21 Th. A]).

\noindent 2. Below, and especially in \S 5, it is helpful if the sets $%
\Sigma _{m}$ are not only in $\mathcal{D}_{\mathcal{B}}$ but also dense. So,
in particular, sets that are locally co-meagre come to mind; however, any
Baire set that is locally co-meagre is co-meagre. (For $\Sigma $ Baire, its
quasi-interior -- the largest (regular) open set equal to $\Sigma $ modulo a
meagre set -- is then locally dense, so everywhere dense and so co-meagre.)

\noindent 3. For $f$ Baire, $f|V$ is continuous in the usual sense (i.e. $%
\mathcal{E\rightarrow E}$) on a $\mathcal{D}_{\mathcal{B}}$-open set $V$
[Oxt, Th. 8.1].

\bigskip

Our approach below is via the \textit{Steinhaus-Weil property} of certain
non-negligible sets $\Sigma $: $0$ is a (usual) interior point of $\Sigma
-\Sigma $. Our motivation comes from some infinite combinatorics going back
to Kestelman [Kes] in 1947 that has later resurfaced in the work of several
authors: Kemperman [Kem] in 1957, Borwein and Ditor [BorD] in 1978, Trautner
[Trau] in 1987, Harry Miller [Mil] in 1989, Grosse-Erdmann [GroE2] in 1989,
and [BinO1,2,3,7] from 2008. The \textit{Kestelman-Borwein-Ditor Theorem}
(KBD below) asserts that \textit{for any Baire/measurable non-negligible }$%
\Sigma $\textit{\ and any null sequence }$z_{n}\rightarrow 0,$\textit{\
there are }$t\in \Sigma $\textit{\ and an infinite} $\mathbb{M}$ \textit{%
such that }$t+z_{m}\in \Sigma $\textit{\ for} $m\in \mathbb{M}$.

On $\mathbb{R}$, KBD\ is both a consequence and a sharpening of the Baire
Category Theorem (BC below). For, BC implies KBD, and conversely -- the
proof of KBD\ requires a sequence of applications of BC [MilO]. The power of
these ideas is shown in the proof of the Uniform Convergence Theorem of
regular variation ([BinGT Ch. 1], [BinO2]).

None of this is special to $\mathbb{R}$: one can work in a Polish abelian
group. Then KBD in this setting implies as an almost immediate consequence
the Effros Theorem ([Ost3], cf. [vMil]), and so the Open Mapping Theorem
[Ost5], as well as other classical results, for instance the
Banach-Steinhaus Theorem -- see the survey [Ost4] and the more recent
developments in [BinO14] and [BinO15, Th. 2].

The significance of the KBD\ is three-fold.

Firstly, if KBD\ applies for the non-negligible sets $\Sigma $ of some
family of sets, then these sets have the Steinhaus-Weil property. For if
not, choose $z_{n}\notin \Sigma -\Sigma $ with $z_{n}\rightarrow 0$
(henceforth termed a `null' sequence); now there are $t\in \Sigma $\textit{\
}and an infinite $\mathbb{M}$ such that $t+z_{m}\in \Sigma $ for $m\in
\mathbb{M}$, so $z_{m}=(t+z_{m})-t\in \Sigma -\Sigma ,$ a contradiction.

Secondly, the several proofs of KBD rely on elementary induction, i.e.
recursion through the natural numbers via DC (see \S 1). As a result our
Berz-type theorems depend only on DC rather than on the full strength of AC
used by Berz.

Finally, any application of KBD in a topological vector space context may be
deemed to take place in the separable subspace generated by the null
sequence.

\bigskip

In an infinite-dimensional separable Banach space: we cannot rely on Haar
measure, as here that does not exist; but we can nevertheless rely on a $%
\sigma $-ideal of sets whose `negligibility' is predicated on the Borel
probability measures of that space\footnote{%
Also relevant here is their regularity: for their outer regularity
(approximation by open sets) see [Par, Th. II.1.2], and their inner
regularity (approximation by compacts) see [Par, Ths. II.3.1 and 3.2], the
latter relying on completeness.}. We recall below their definition and two
key properties, the first of which relies on \textit{separability} (hence
the frequent recourse below to separable Banach subspaces): the
Steinhaus-Weil property and the \textit{weak extension of the Fubini theorem}
(WFT; see below)\ due to Christensen [Chr1], which may be applied here. For
this we need to recall that $B\subseteq G$ is \textit{universally measurable}
if $B$ is measurable with respect to every Borel measure on $G$ -- for
background, see e.g. cf. [Fre2, 434D, 432]. Examples are analytic subsets
(see e.g. [Rog, Part 1 \S 2.9], or [Kec, Th. 21.10], [Fre2, 434Dc]) and the $%
\sigma $-algebra that they generate. Beyond these are the \textit{provably }$%
\mathbf{\Delta }_{2}^{1}$ sets of [FenN], defined in Appendix 1.1 below.

The $\sigma $-ideal of \textit{Haar null} sets is a generalization of
Christensen [Chr1,2] to a non-locally compact group of the notion of a Haar
measure-zero set: see again Hoffmann-J\o rgensen [Rog, Part 3, Th. 2.4.5]
and Solecki [Sole1,2,3] (and [HunSY] in the function space setting).

Christensen [Chr1] shows that in an abelian Polish group $(G,\cdot )$ the
family $\mathcal{H}(G,\cdot )$ of Haar null sets forms a $\sigma $-ideal.
This was extended for \textit{`}left Haar null' sets (see below) by Solecki
[Sole3, Th. 1] in the more general setting of (not necessarily abelian)
Polish groups $(G,\cdot )$ \textit{amenable at 1}, the scope of which he
studies, in particular proving that any abelian Polish group is amenable at
1 [Sole3, Prop 3.3]; this includes, as additive groups, separable $F$- (and
hence Banach) spaces.

A subset of a Polish group $G$ is \textit{left Haar null} [Sole3] if it is
contained in a universally measurable set $B$ (for which see \S 3) such that
for some Borel probability measure $\mu $ on $G$%
\[
\mu (gB)=0\qquad (g\in G).
\]%
Solecki also considers the (in general) narrower family of Haar null sets
(as above, but now $\mu (gBh)=0$ for all $g,h\in G).$ Below we work in
vector spaces and so the non-abelian distinctions vanish.

The Steinhaus-Weil Theorem holds also for universally measurable sets that
are not Haar null; this was proved by Solecki (actually for left Haar null
[Sole3, Th. 1(ii)]; cf. Hoffmann-J\o rgensen [Rog, Part 3, Th. 2.4.6]) by
implicitly proving KBD. It may be checked that his proof uses only DC. One
may also show that the KBD\ theorem follows from amenability at 1: see
[BinO15].

Christensen's WFT [Chr1] (for a detailed proof, see [BorM]) concerns the
product $H\times T$ of a locally compact group $T,\ $equipped with Haar
measure $\eta ,$ and an arbitrary abelian Polish group $H,$ and is a
`one-way round' theorem (for $T$-sections): if $A\subseteq H\times T$ is
universally measurable, then $A$ is Haar null iff the sections $A(h):=\{t\in
T:(h,t)\in A\}$ are Haar measure-zero except possibly for a Haar null set
(in the sense above) of $h\in H$ (i.e. for `quasi all' $h\in H$). The `other
way round' (for $H$-sections$)$ may fail, as was shown by Christensen [Chr1,
Th. 6].

\section{Sublinearity and Berz's Theorem}

We begin with Theorem 2 on subadditive functions. We deal with the category
and measure versions together via the Steinhaus-Weil property, and use DC
rather than AC. The Lebesgue measurable case (with AC) is classical [HilP,
Ch. 7]. See also [BinO1, Prop. 1], [Kuc, Th. 16.2.2], [BinO13, Prop. 7$%
^{\prime }$].

Recalling that $\mathcal{T}$ is a Baire space topology if Baire's theorem
holds under $\mathcal{T}$, say that a vector space $X$ has a \textit{%
Steinhaus-Weil topology }$\mathcal{T}$ if the \textit{non-meagre Baire sets}
of $\mathcal{T}$ \textit{have the Steinhaus-Weil property}. Thus $\mathcal{E}
$ and $\mathcal{D}_{\mathcal{L}}$ are Steinhaus-Weil topologies for $\mathbb{%
R}$ that are Baire topologies, by the classical Steinhaus-Piccard-Pettis
Theorems (see e.g. [BinO9]).

As above, say that $S:X\rightarrow \mathbb{R}$ is $\mathcal{T}$-Baire if $%
S^{-1}$ takes (Euclidean) open sets of $\mathbb{R}$ to $\mathcal{T}$-Baire
sets in $X$. Thus $\mathcal{E}$-Baire means Baire in the usual sense, and $%
\mathcal{D}_{\mathcal{L}}$-Baire means Lebesgue measurable.

\bigskip

\noindent \textbf{Theorem 2}. \textit{For }$X$ \textit{a vector space with a
Steinhaus-Weil, Baire topology }$\mathcal{T}$\textit{\ and} $S:X\rightarrow
\mathbb{R}$\textit{\ subadditive:} \textit{if} $S$ \textit{is }$\mathcal{T}$%
\textit{-Baire, then it is locally bounded.}

\bigskip

\noindent \textbf{Proof.} Suppose $|S(u+z_{n})|\rightarrow \infty $ for some
$u\in X$ and null sequence $z_{n}\rightarrow 0.$ As the level sets $%
H_{n}^{\pm }:=\{t:|S(\pm t)|\leq n\}$ are $\mathcal{T}$-Baire and $\mathcal{T%
}$ is Baire, for some $k$ the set $H_{k}^{\pm }$ is non-meagre. As $\mathcal{%
T}$ is Steinhaus-Weil, $H_{k}^{\pm }-H_{k}^{\pm }$ has 0 in its interior. So
there is $n\in \mathbb{N}$ such that $z_{m}\in H_{k}^{\pm }-H_{k}^{\pm }$
for all $m\geq n.$ For $m\geq n,$ choose $a_{m},b_{m}\in H_{k}^{\pm }$ with
\[
z_{m}=a_{m}-b_{m}
\]%
for $m\geq n.$ Then for all $m\geq n$%
\begin{eqnarray*}
S(u)-2k &\leq &S(u)-S(-a_{m})-S(b_{m})\leq S(u+a_{m}-b_{m}) \\
&=&S(u+z_{m})\leq S(u)+S(a_{m})+S(-b_{m})\leq S(u)+2k,
\end{eqnarray*}%
contradicting unboundedness. $\square $

\bigskip

The key to Theorem 1B is Theorem 3 below. It may be regarded as a
subadditive analogue of Ostrowski's Theorem for additive functions (cf.
[BinO9], [MatS]). The result extends its counterpart in [BinO10, Prop. 13],
with a simpler proof, and uses only DC (via KBD). As it depends on the
Steinhaus-Weil property, it handles the Baire and measurable cases together.
In the theorem below, we write $\mathbb{R}_{>}:=\mathbb{R}_{+}\backslash
\{0\}$ and similarly $\mathbb{R}_{<}$.

\bigskip

\noindent \textbf{Theorem 3.} \textit{If }$S:\mathbb{R}\rightarrow \mathbb{R}
$\textit{\ is subadditive, locally bounded with }$S(0)=0,$ \textit{and:}

\noindent (i)\textit{\ there is a symmetric set }$\Sigma $\textit{\ (i.e. }$%
\Sigma =-\Sigma $\textit{) containing }$0$ \textit{with }$S|\Sigma $\textit{%
\ continuous at }$0$\textit{;}

\noindent (ii) \textit{for each }$\delta >0,$ $\Sigma _{\delta }^{+}:=\Sigma
\cap (0,\delta )$ \textit{has the Steinhaus-Weil property}\newline
-- \textit{then }$S$ \textit{is continuous at }$0$\textit{\ and so
everywhere. }

\textit{In particular, this is so if }$S(0)=0$ \textit{and there is a
symmetric set }$\Sigma $\ \textit{containing }$0$\textit{\ on which}%
\[
S(u)=c^{\pm }u\text{ }\mathit{for\ some\ }c^{\pm }\in \mathbb{R}\mathit{\
and\ all\ }u\in \Sigma \cap \mathbb{R}_{>}\mathit{,or}\text{ }\mathit{all}%
\text{ }u\in \Sigma \cap \mathbb{R}_{<}\text{ }\mathit{resp.,}
\]%
\textit{and }$\Sigma $ \textit{is} \textit{Baire/measurable, non-negligible
in each }$(0,\delta )$\textit{\ for }$\delta >0$\textit{.}

\bigskip

\noindent \textbf{Proof.} If $S$ is not continuous at $0$, then (see e.g.
[HilP, 7.4.3], cf. [BinO13, Prop. 7]) $\lambda _{+}:=\lim \sup_{t\rightarrow
0}S(t)>\lim \inf_{t\rightarrow 0}S(t)\geq 0,$ the last inequality by
subadditivity and local boundedness at 0. Choose $z_{n}\rightarrow 0$ with $%
S(z_{n})\rightarrow \lambda _{+}>0.$ Let $\varepsilon =\lambda _{+}/3.$ By
continuity on $\Sigma $ at $0,$ there is $\delta >0$ with $%
|S(t)|<\varepsilon $ for $t\in \Sigma \cap (-\delta ,\delta ).$ By the
Steinhaus-Weil property of $\Sigma _{\delta }$ there is $n$ such that $%
z_{m}\in \Sigma _{\delta }^{+}-\Sigma _{\delta }^{+}$ for all $m\geq n.$
Choose $a_{m},b_{m}\in \Sigma _{\delta }^{+}$ with $z_{m}=a_{m}-b_{m};$ so
by subadditivity
\[
S(z_{m})\leq S(a_{m})+S(-b_{m})\leq 2\varepsilon .
\]%
Taking limits,%
\[
\lambda _{+}\leq 2\varepsilon <\lambda _{+}.
\]%
This contradiction shows that $S$ is continuous at $0.$ As in [HilP, Th.
2.5.2], continuity at all points follows by noting that%
\[
S(x)-S(-h)\leq S(x+h)\leq S(x)+S(h).
\]

The remaining assertion follows from the Steinhaus-Piccard-Pettis Theorem
via Theorem 2, as continuity at $0$ on $\Sigma $ follows from%
\[
|c^{\pm }u|\leq |u|\cdot \max \{|c^{+}|,|c^{-}|\}.\qquad \square
\]

\bigskip

\noindent \textbf{Proofs of Theorems 1B and 1M. }Let $S:\mathbb{R}%
\rightarrow \mathbb{R}$ be sublinear and either Baire or measurable, that
is, Baire in one of the two topologies $\mathcal{D}_{\mathcal{B}}(\mathcal{E}%
)$ or $\mathcal{D}_{\mathcal{B}}(\mathcal{D}_{\mathcal{L}}).$ By Theorem BL $%
S$ is quasi $\sigma $-continuous. Taking $\Sigma _{m}$ as in the Definition
in \S 2 (with $m$ fixed), apply Lemma S to $A=B=\Sigma _{m}\subseteq \mathbb{%
R}_{+}.$ Fix (non-zero!) $a,b\in \Sigma _{m};$ as these are density points,
there are $a_{n},b_{n}\in \Sigma _{m}$ and $q_{n}\in \mathbb{Q}_{+}$ so that
\[
b_{n}=q_{n}a_{n};\quad a=\lim\nolimits_{n}a_{n},\text{ }b=\lim%
\nolimits_{n}b_{n}.
\]%
As $S$ is sublinear,
\[
S(q_{n}a_{n})/S(a_{n})=q_{n}=b_{n}/a_{n}\rightarrow b/a.
\]%
But $S|\Sigma _{m}$ is continuous at $a$ and $b,$ so%
\[
S(b)/S(a)=\lim\nolimits_{n}S(b_{n})/S(a_{n})=b/a.
\]%
So on $\Sigma _{m},$ $S(u)=c_{m}u.$ But $\Sigma _{m}\supseteq \Sigma _{0},$
so $c_{m}=c_{0}$ for all $m.$ So $S(x)=c_{0}x$ for $x\in \Sigma ^{+},$ i.e.
for almost all $x>0.$ Repeat for $\mathbb{R}_{-}$ with an analogous set $%
\Sigma ^{-}$ and put $\Sigma :=\{0\}\cup \Sigma ^{+}\cup \Sigma ^{-}.$ We
may assume $-\Sigma =\Sigma $ (otherwise pass to the subset $\Sigma \cap
(-\Sigma ),$ which is quasi all of $\mathbb{R}$). By Th. 3, $S$ is
continuous at $0$ and so everywhere. In summary: $S$ is linear on the
\textit{dense} subset $\Sigma ^{+}\subseteq \mathbb{R}_{+}$ and continuous,
and likewise on the \textit{dense} subset $\Sigma ^{-}\subseteq \mathbb{R}%
_{-}.$ So $S$ is linear on the whole of $\mathbb{R}_{+}$, and on the whole
of $\mathbb{R}_{-}$. $\square $

\bigskip

For later use (in \S 4 below) we close this section with a variant on
Theorem 2, Theorem 2H (`H for Haar'). The proof is the same as that of
Theorem 2 above, but needs a little introduction. Say that a function $%
S:G\rightarrow \mathbb{R}$ is \textit{universally measurable} if $S^{-1}$
takes open sets to universally measurable sets (as in \S 2) in $G$; further
say that a $\sigma $-ideal $\mathcal{H}$ of subsets of a topological vector
space $X$ (the `negligible sets') is \textit{proper} if $X\notin \mathcal{H}$%
, and that $\mathcal{H}$ has the \textit{Steinhaus-Weil property} if
universally measurable sets that are not in $\mathcal{H}$ have the interior
point property.

\bigskip

\noindent \textbf{Theorem 2H}. \textit{For }$X$ \textit{a topological vector
space, }$\mathcal{H}$\textit{\ a proper }$\sigma $\textit{-ideal with the
Steinhaus-Weil property\ and} $S:X\rightarrow \mathbb{R}$\textit{\
subadditive:} \textit{if} $S$ \textit{is universally measurable, then it is
locally bounded.}

\section{Banach versions}

In Theorem 4B and 4M below we extend the category and measure results in
Theorems 1B and 1M to the setting of a Banach space $X.$ Since the
conclusions are derived from continuity (and local boundedness), our results
are first established for separable (sub-) spaces, which then extend to the
non-separable context, by Theorem B (\S 1). The key in each case is an
appropriate application of Theorem 2 (or Theorem 2H). The category case here
is covered by the Piccard-Pettis Theorem, true for non-meagre Baire sets in $%
X$; in fact more is true, as KBD\ holds in any analytic group with a
translation-invariant metric -- see [BinO8, Ths. 1.2, 5.1] or [Ost2, Th. 2],
which also covers $F$-spaces, so including Fr\'{e}chet spaces (see the end
of this section). In the absence of Haar measure, the analogous `measure
case' arising from universal measurability is technically more intricate,
but nevertheless true -- see below. It is here that our methodology requires
us to \textit{pass down} to separable subspaces of a Banach space $X$. That
this suffices to reduce the case of a general Banach space $X$ to the
separable case follows from the result below, a corollary of Theorem B of \S %
1. Henceforth we write $B_{\delta },\Sigma _{\delta }$ respectively for the
closed unit ball $\{x:||x||\leq \delta \}$ and the $\delta $-sphere $%
\{x:||x||=\delta \},$ and use the following notation for \textit{lines} and
\textit{rays}:%
\[
R(u):=\{\lambda u:\lambda \in \mathbb{R}\},\quad R_{\pm }(u):=\{\lambda
u:\lambda \in \mathbb{R}_{\pm }\cup \{0\}\}.
\]

\bigskip

\noindent \textbf{Corollary B.} \textit{For }$X$\textit{\ a Banach space and}
$S:X\rightarrow \mathbb{R}$ \textit{subadditive, if }$S|B$\textit{\ is
locally bounded for each closed separable subspace }$B$\textit{, then }$%
\{|S(x)|/||x||:x\neq 0\}$\textit{\ is bounded.}

\bigskip

\noindent \textbf{Proof. }By Theorem B, $S$ is continuous on $X$ , so there
is $\delta >0$ with%
\[
||S(v)||\leq 1\qquad (||v||\leq \delta ).
\]%
Furthermore, for any $x\neq 0$ taking $u:=x/||x||,$ the restriction of $S$
to the ray $R_{+}(u)$ is positively homogeneous by Theorem 1B, and so%
\[
|S(x)|=|S(||x||u)|=|S(\delta u)||x||/\delta |\leq ||x||/\delta .\qquad
\square
\]

\bigskip

The proof of the category case in Theorem 4B below would have been easier
had we used AC to construct the function $c(x);$ but, as we wish to rely
only on DC, more care is needed.

We offer two proofs. The first uses Theorems 1B and 2 (and is laid out so as
to extend easily to the more demanding $F$-space setting of Theorem 4F\
below); the second is more direct, but uses a classical selection
(uniformization) theorem, together with a Fubini-type theorem for negligible
sets in a product space. Both proofs have Banach-space `measure' analogues.

\bigskip

\noindent \textbf{Theorem 4B.} \textit{For }$X$\textit{\ a Banach space, and
}$S:X\rightarrow \mathbb{R}$\textit{\ Baire, if }$S$\textit{\ is subadditive
and }$\mathbb{Q}_{+}$\textit{-homogeneous, then}

\noindent\ (i)\textit{\ }$S$\textit{\ is continuous and convex with epigraph
a convex cone pointed at }$0$\textit{, and }

\noindent\ (ii) \textit{there is a bounded function }$c:X\rightarrow \mathbb{%
R}$\textit{\ such that}%
\[
S(x)=c(x)||x||.
\]

\noindent \textit{First Proof. }Since $S$ is mid-point convex, and we first
seek to establish continuity, we begin by establishing it for any separable
subspace; we then use Theorem and Corollary B above to draw the same
conclusion about $X$ itself. Consequently, we may w.l.o.g. assume $X$ is
itself separable. By Theorem 2 applied to the usual meagre sets, $S$ is
locally bounded at $0,$ so there are $M$ and $\delta >0$ such that%
\[
|S(x)|\leq M\qquad (x\in B_{\delta }).
\]%
In particular, for $v\in \Sigma _{\delta },$ $|S(v)|\leq M.$ For $u\in X$
define a ray-restriction of $S$ by%
\[
f_{u}(x):=S(x)\quad (x\in R(u)).
\]%
For fixed $u,$ as the mapping $(\lambda ,u)\longmapsto \lambda u$ from $%
\mathbb{R}$ into $X$ is continuous, the set $R(u)$ is $\sigma $-compact. So
for any fixed $u,$ $f_{u}$ is Baire; indeed, $f_{u}(x)\in (a,b)$ iff $%
S(x)\in (a,b)$ and $x\in R(u),$ i.e.%
\[
\{x:f_{u}(x)\in (a,b)\}=\{x:S(x)\in (a,b)\}\cap R(u),
\]%
and $R(u)$ has the Baire property (being $\sigma $-compact). So by Th. 1B,
for any fixed $u$ the function $f_{u}$ is continuous and there exist $c^{\pm
}\in \mathbb{R}$ with $S(\lambda u)=c^{\pm }\lambda $ according to the sign
of $\lambda .$ This justifies the definitions below for $u\in X$:%
\[
c^{+}(u):=S(u),\qquad c^{-}(u):=-S(-u).
\]%
Then, for fixed $u,$ by continuity of $f_{u},$ $S(\lambda u)=c^{+}(u)\lambda
=\lambda S(u)$ for $\lambda \geq 0,$ so that $S$ is \textit{positively
homogeneous} on $R_{+}(u);$ likewise, $S(\lambda u)=S((-\lambda
)(-u))=(-\lambda )S(-u)=c^{-}(u)\lambda $ for $\lambda \leq 0.$ Then for $%
u=x/||x||$ with $x\neq 0,$ as $v:=\delta u\in \Sigma _{\delta },$
\[
|S(x)|=|S(||x||u)|=|S(v)|\cdot ||x||/\delta \leq (M/\delta )||x||.
\]%
So $S$ is continuous at $0,$ and so by subadditivity everywhere, as in the
proof of Theorem 3. By continuity (as $S$ is positively homogeneous) $S$ is (%
$\mathbb{R}$-) convex [Roc, Th. 4.7]; so its epigraph is a convex cone
pointed at the origin [Roc, Th. 13.2].

Finally, for $x\neq 0,$ take $c(x):=S(x/||x||),$ which as above is bounded
by $M/\delta ;$ then
\[
S(x)=c(x)||x||.\qquad \square
\]

\bigskip

\noindent \textit{Second Proof. }As above, we again assume that $X$ is
separable. By Theorem BL (\S 2) there is a co-meagre subset $\Sigma $ with $%
S|\Sigma $ continuous. By passage to $\Sigma \cap (-\Sigma )$ we may assume $%
\Sigma $ is symmetric; we may also assume that $\Sigma $ is a $\mathcal{G}%
_{\delta }$. By the Kuratowski-Ulam Theorem [Oxt, Th. 15.1], for quasi all $%
x\neq 0,$ say for $x\in D$ with $D$ a $\mathcal{G}_{\delta }$-set, the ray $%
R_{+}(x)\cap \Sigma $ is co-meagre on $\Sigma .$ By the Steinhaus Theorem,
Sierpi\'{n}ski's Lemma S applies. By Theorem 1B for $s\in \Sigma \cap
R_{+}(x)$ there is $c$ with $S(s)=c||s||$ (as $s=x||s||/||x||).$ Now $%
S|\Sigma $ is continuous so a Borel function, as $\Sigma $ is a $\mathcal{G}%
_{\delta }$, and for fixed $x,$ $S(a)/||a||$ is constant for (density)
points $a$ of $\Sigma \cap R_{+}(x)$ for $x\in D.$ (This uses the isometry
of $R_{+}(x)$ and $\mathbb{R}_{+}.)$ So by Novikov's Theorem (see e.g.
[JayR], p. x] -- cf. [Kec, 36.14]) there is a Borel function $c:D\rightarrow
\mathbb{R}$ such that $S(x)=c(x)||x||$ for $x\in D.$ By Theorem 1B, since $S$
is bounded near the origin, $c(x)$ is also bounded on $D$ near $0$ (as in
the previous proof). From this boundedness near $0,$ by Theorem 3, $S(x)$ is
continuous for all $x$. By continuity, $S$ is positively homogeneous, so
again convex with epigraph a convex cone pointed at the origin. $\square $

\bigskip

\noindent \textbf{Remarks.} 1. In the first proof, one may show that $S$ is
continuous at $0$ by considering a (null) non-vanishing sequence $%
z_{n}\rightarrow 0.$ Put $u_{n}:=z_{n}/||z_{n}||;$ by DC select $c_{n}^{\pm
} $ such that $S(\lambda u_{n})=c_{n}^{\pm }\lambda $, according to the sign
of $\lambda $. As $S$ is locally bounded at $0,$ there are $M$ and $\delta
>0 $ such that%
\[
|S(x)|\leq M\qquad (x\in B_{\delta }).
\]%
W.l.o.g. $\delta \in \mathbb{Q}_{+}$, so for $x=\delta u_{n}\in B_{\delta },$
$|S(x)|=|c_{n}^{+}\delta |\leq M;$ then $|c_{n}^{+}|\leq M/\delta .$ So%
\[
S(z_{n})=S(||z_{n}||u_{n})=c_{n}^{+}||z_{n}||\leq (M/\delta
)||z_{n}||\rightarrow 0.
\]

\noindent 2. The second proof uses the Fubini-like Kuratowski-Ulam Theorem
[Oxt, 15.1] (cf. [Chr1]). This can fail in a non-separable metric context,
as shown in [Pol] (cf. [vMilP]), but see [FreNR] and [Sole4].

\bigskip

Either argument for Theorem 4B above has an immediate \textit{Lebesgue
measure} analogue for $X=\mathbb{R}^{d},$ and beyond that a \textit{Haar
measure} analogue for $X$ a locally compact group with Haar measure $\eta ,$
by the classical Fubini Theorem (see e.g. [Oxt, Th. 14.2]). But we may reach
out further still for a measure analogue, Theorem 4M\ below, by employing
the $\sigma $-ideal of Haar null sets (\S 2). Whilst our argument is simpler
(through not involving radial open-ness), there is a close relation to the
result of [FisS], which is concerned with \textit{convex} functions $S$ that
are measurable in the following sense: $S^{-1}$ takes open sets to sets
that, modulo Haar null sets, are universally measurable sets in $X$ (we turn
to convexity in \S 5: see especially Th. 7 and 8). Below (recall \S 3) a
function $S:G\rightarrow \mathbb{R}$ is \textit{universally measurable} if $%
S^{-1}$ takes open sets to universally measurable sets in $G$; this means
that, as in \S 3, the level sets $H_{n}^{\pm }$ are universally measurable,
so if $G\ $is amenable at 1, in particular if $G$ is an abelian Polish
group, for some $k\in \mathbb{N}$ the level set $H_{k}^{\pm }$ is not Haar
null (since $X=\bigcup\nolimits_{n\in \mathbb{N}}H_{n}^{\pm }$ is not Haar
null). This aspect would remain unchanged if the level sets were universally
measurable modulo Haar null sets.

\bigskip

\noindent \textbf{Theorem 4M. }\textit{For }$X$\textit{\ a Banach space, and
}$S:X\rightarrow \mathbb{R}$\textit{\ universally measurable: if }$S$\textit{%
\ is subadditive and }$\mathbb{Q}_{+}$\textit{-homogeneous, then}

\noindent\ (i)\textit{\ }$S$\textit{\ is continuous and convex with epigraph
a convex cone pointed at }$0$\textit{, and }

\noindent\ (ii) \textit{there is a bounded function }$c:X\rightarrow \mathbb{%
R}$\textit{\ such that}%
\[
S(x)=c(x)||x||.
\]

\noindent \textit{First Proof. }Proceed as in the first proof of Theorem 4B
(reducing as there to separability), but in lieu of Theorem 2 apply Theorem
2H here to the $\sigma $-ideal of Haar null sets $\mathcal{H}(X,+)$. $%
\square $

\bigskip

\noindent \textit{Second Proof. }Reduce as before to the separable case.
With WFT\ above, as a replacement for the Kuratowski-Ulam theorem, we may
follow the proof strategy in the second proof of Theorem 4B, largely
verbatim. Regarding the line $R(x)$ (for $x\neq 0$) as a locally compact
group isomorphic to $\mathbb{R}$, take $\mu :=\mu _{\Sigma }\times \eta _{1}$
to be a probability measure with atomless spherical component $\mu _{\Sigma
} $ (a probability on the unit sphere of $X;$ this can be done since the
atomless measures form a dense $\mathcal{G}_{\delta }$ under the weak
topology -- cf. [Par, Th. 8.1]) and radial component $\eta _{1}$ a \textit{%
probability} on $\mathbb{R}$ absolutely continuous with respect to Lebesgue
(Haar) measure. We claim that $S|R_{+}(x)$ is quasi-$\sigma $-continuous on
a (Haar/Lebesgue) co-null set for quasi all $x.$ For if not, there is a set $%
C$ that is not Haar null with $S|R_{+}(x)$ not $\sigma $-continuous for $%
x\in C$. So there is $u\in X$ with $\mu (u+C)>0,$ and so $u+C\ $is not
radial. Put $m(B):=\mu (u+B)$ for Borel sets $B,$ again a probability
measure. By Theorem BL and WFT, $S|R_{+}(x)$ is quasi-$\sigma $-continuous
for $m$-almost all $x$, except on some set $E$ with $m(E)=0.$ This is a
contradiction for points in $C\backslash E.$ Now continue as in Theorem 4B. $%
\square $

\bigskip

\noindent \textit{F-spaces. }Recall that an $F$\textit{-space} is a
topological vector space with topology generated by a complete
translation-invariant metric $d_{X}$ ([KalPR, Ch. 1], [Rud, Ch. 1]). Thus
the topology is generated by the $F$\textit{-norm} $||x||:=d_{X}(0,x),$
which satisfies the triangle inequality with $||\alpha x||\leq ||x||$ for $%
|\alpha |\leq 1,$ and under it scalar multiplication is jointly continuous.
This continuity implies that a vector $x$ can be scaled down to arbitrarily
small size. Consequently, the proofs above may be re-worked to yield F-space
versions\textit{\ }of Theorems 4B and 4M. However, in the absence of a norm
(see \S 6.5 for normability), there is no isometry between the rays $R(x)$
below and $\mathbb{R}_{+}$, only an injection $\Delta :R(x)\rightarrow
\mathbb{R}_{+}$. We are thus left with a result that has a somewhat weaker
representation of $S$.

We need the $F$-norm to be \textit{unstarlike}, in the sense that the
`norm-length' (i.e. the range of the norm) of all rays be the same, say
unbounded for convenience. This last property holds for the $L^{p}$ spaces
for $0<p<1$ with the familiar $F$-norm $||f||:=\left( \int
|f(t)|^{p}dt\right) ^{1/p}.$

Unstarlikeness is an $F$-norm, rather than a topological, property; it will
hold after re-norming, albeit with $(0,1)$ as the common range$,$ when
taking the $F$-norm to be $||x||:=\sup_{n}2^{-n}(||x||_{n}/(1+||x||_{n})),$
for $||\cdot ||_{n}$ a \textit{distinguishing} sequence of semi-norms, since
$\varphi _{x,n}(t):=t||x||_{n}/(1+t||x||_{n})$ maps $\mathbb{[}0,\infty )$
onto $[0,1).$ Examples here are provided by spaces of continuous functions
such as $C(\Omega ),$ for $\Omega :=\bigcup_{n}K_{n}$ with $K_{n}\subseteq $%
\textrm{int}($K_{n+1})$ a chain of compact subsets of a Euclidean space, and
with $||f||_{n}:=||f|_{K_{n}}||_{\infty }$ for $||\cdot ||_{\infty }$ the
supremum norm. Likewise this holds in the subspace $H(\Omega )$ of
holomorphic functions, and in $C^{\infty }(\Omega )$ when $||f||_{n}:=\max
\{||D^{\alpha }f||_{\infty }:|\alpha |<n\}$ for multi-indices $\alpha $ --
see [Rud, \S 1.44-47]). Being infinite-dimensional, none of them are
normable as they are either locally bounded or Heine-Borel (or both) -- cf.
[Rud, Th. 1.23].

\bigskip

\noindent \textbf{Theorem 4F.} \textit{For }$X$\textit{\ an }$F$\textit{%
-space and }$S:X\rightarrow \mathbb{R}$\textit{\ Baire, if }$S$\textit{\ is
subadditive and }$\mathbb{Q}_{+}$\textit{-homogeneous, then}

\noindent\ (i)\textit{\ }$S$\textit{\ is continuous and convex with epigraph
a convex cone pointed at }$0$\textit{, and }

\noindent\ (ii) \textit{for any unstarlike }$F$\textit{-norm }$||\cdot ||$%
\textit{\ (with }$||tx||\rightarrow \infty $ $(t\rightarrow \infty )$
\textit{for all }$x\neq 0),$ \textit{there are a bounded function }$%
c:X\rightarrow \mathbb{R}$\textit{, a constant }$\delta ,$ \textit{and an
injection }$\Delta :R(x)\rightarrow \mathbb{R}_{+}$\textit{\ such that}%
\[
S(x)=c(x)\Delta (x),\text{ where }||x/\Delta (x)||=\delta .
\]%
\textit{In particular, if }$X$\textit{\ is normable with norm }$||\cdot
||_{X}$, \textit{then }$\Delta (x)=||x||_{X}/\delta .$

\bigskip

\noindent \textit{Proof. }Let $||.||$ be an unstarlike $F$-norm. For any $%
x\neq 0,$ the map $\varphi _{x}:t\mapsto ||tx||$ is a continuous injection
with $\varphi _{x}(0)=0$ and $\varphi _{x}(1)=||x||$; so for $||x||\geq
\delta $ we may define $\delta (x):=\min \{t:||tx||=\delta \},$ the infimum
being attained. So $||\delta (x)x||=\delta .$ The unstarlike property
implies that $\delta (x)$ is likewise well defined for all $x\neq 0.$

We now assume w.l.o.g. that $X$ is separable, as in the earlier variants of
Th. 4, for the same reasons (though we need the $F$-norm analogue of
Corollary B, also valid -- see below for the relevant positive homogeneity).
Proceed as in the first proofs of Theorems 4B and 4M, with a few changes,
which we now indicate. Of course we refer respectively to the $\sigma $%
-ideals of meagre sets and Haar null sets.

With this in mind one deduces again positive homogeneity, and thence, for $%
x\neq 0$ and with $v=\delta (x)x\in \Sigma _{\delta },$ that as $\delta
(x)>0 $
\[
|S(x)|=|S(v/\delta (x))|=|S(v)|/\delta (x)\leq M/\delta (x).
\]%
Now $\delta (x)\rightarrow \infty $ as $x\rightarrow 0,$ and so $S$ is
continuous at $0;$ indeed, for each $n\in \mathbb{N}$ the function $x\mapsto
||nx||$ is continuous at $x=0,$ so by DC there is a positive sequence $%
\{\eta (n)\}_{n\in \mathbb{N}}$ such that $||nx||<\delta $ for all $x\in
B_{\eta (n)}.$ So $\delta (x)>n$ for $x\in B_{\eta (n)}$ and $n\in \mathbb{N}%
,$ and so%
\[
|S(x)|\leq M/\delta (x)<M/n\qquad (x\in B_{\eta (n)}).
\]%
Thereafter, taking $c(x):=c^{+}(\delta (x)x)=S(\delta (x)x),$ which is
bounded by $M,$
\[
S(x)=S(\delta (x)x/\delta (x))=c(x)\Delta (x),\text{ where }\Delta
(x):=1/\delta (x).
\]%
So $||x/\Delta (x)||=\delta .$ If the $F$-norm is a norm, $\delta
(x):=\delta /||x||;$ then $||x\delta (x)||=\delta ,$ so that $\Delta
(x):=||x||/\delta .$ $\square $

\bigskip

Theorem 4F implies Theorem 4B and 4M by taking $c(x)/\delta $ in place of $%
c(x).$

\section{Convexity}

We begin by recalling a classical result, Theorem BD below, which motivates
the themes of this section. These focus on the two properties of a function $%
S$ of \textit{mid-point convexity}%
\[
S\left( \frac{1}{2}(x+y)\right) \leq \frac{1}{2}\left( S(x)+S(y)\right) ,
\]%
and \textit{convexity, }which, for purposes of emphasis, we also refer to
(as in [Meh]) as \textit{full }(or $\mathbb{R}$-) \textit{convexity}:%
\[
S((1-t)x+ty)\leq (1-t)S(y)+tS(y)\qquad (t\in (0,1)),
\]%
by considering the weaker property of \textit{mid-point convexity on a set} $%
\Sigma :$%
\[
S\left( \frac{1}{2}(x+y)\right) \leq \frac{1}{2}\left( S(x)+S(y)\right)
\qquad (x,y\in \Sigma ).
\]%
This is the \textit{weak mid-point convexity} of \S 1.

In Theorems 5-7 below, we give \textit{local} results, with the hypotheses
holding on a set $\Sigma .$ The smaller $\Sigma $ is, the more powerful (and
novel) the conclusions are. For instance, $\Sigma $ might be locally
co-meagre (and so co-meagre, as noted in the remarks to the definition of
quasi-$\sigma $-continuity in \S 2).

\bigskip

\noindent \textbf{Theorem BD (Bernstein-Doetsch Theorem, }[Kuc, \S\ 6.4]%
\textbf{).} \textit{For }$X$\textit{\ a normed vector space, if }$%
S:X\rightarrow \mathbb{R}$\textit{\ is mid-point convex and locally bounded
somewhere (equivalently everywhere), then }$S$\textit{\ is continuous and
fully convex.}

\bigskip

\noindent \textbf{Proof. }This is immediate from Theorem B (see \S 1). See
also the `third proof' in [Kuc, \S\ 6.4], as the other two apply only in $%
\mathbb{R}^{d}.$ $\square $

\bigskip

The theorem gives rise to a sharp dichotomy for mid-point convex functions,
similar to that for additive functions: they are either continuous
everywhere or discontinuous everywhere (`totally discontinuous'), since
local boundedness is `transferable' between points. So on the one hand, a
Hamel basis yields discontinuous additive examples (the `Hamel pathology' of
[BinGT, \S 1.1.4]) and, on the other, a smidgen's worth of regularity
prevents this -- see Corollary 1 below -- and the mid-point convex functions
are then continuous.

A closely related result (for which see e.g. [Sim, Prop. 1.18]) we give as
Theorem BD$^{\text{*}}$ below, whose proof we include, as it is so simple.

\bigskip

\noindent \textbf{Theorem BD}$^{\text{*}}$\textbf{.} \textit{For }$X$\textit{%
\ a normed vector space, if }$S:X\rightarrow \mathbb{R}$\textit{\ is fully
convex and locally bounded, then }$S$\textit{\ is continuous.}

\bigskip

\noindent \textbf{Proof. }W.l.o.g. assume that $S$ is bounded in the unit
ball, by $K\ $say (otherwise translate to the origin and rescale the norm).
For $x$ in the unit ball, setting $u=x/||x||,$ and first writing $x$ as a
convex combination of $0$ and $u,$ then $0$ as a convex combination of $-u$
and $x,$%
\[
S(x)<(1-||x||)S(0)+||x||S(u),\qquad (1+||x||)S(0)<||x||S(-u)+S(x).
\]%
From here%
\[
||x||[S(0)-S(-u)]<S(x)-S(0)<||x||[S(u)-S(0)]:\quad |S(x)-S(0)|<2K||x||.\quad
\square
\]

Thus the emphasis in convexity theory is on generic differentiability; for
background see again [Sim]. In Theorem 6 below we derive continuity and full
(i.e. $\mathbb{R}$-) convexity, as in Theorem BD [Kuc, \S\ 6.4], for
functions possessing the weaker property of mid-point convexity on certain
subsets $\Sigma $ of their domain with negligible complement, for instance
co-meagre or co-null sets. The results below vary their contexts between $%
\mathbb{R}$ and a general Banach space, and refer to sets with the following
Steinhaus-Weil property.

\bigskip

\textbf{Definition. }Say that $\Sigma $ is \textit{locally Steinhaus-Weil},
or has the \textit{Steinhaus-Weil property locally}, if for $x,y\in \Sigma $
and, for all $\delta >0$ sufficiently small, the sets $\Sigma
_{z}^{+}:=\Sigma \cap B_{\delta }(z),$ for $z=x,y,$ have the interior point
property that $\Sigma _{x}^{+}-\Sigma _{y}^{+}$ has $x-y$ in its interior.
(Here $B_{\delta }(x)$ is the closed ball about $x$ of radius $\delta .)$
See [BinO7] for conditions under which this property is implied by the
interior point property of the sets $\Sigma _{x}^{+}-\Sigma _{x}^{+}$ (cf.
[BarFN]).

\bigskip

Examples of locally Steinhaus-Weil sets relevant here are the following:

\noindent (i) $\Sigma $ density-open in the case $X:=\mathbb{R}^{n}$ (by
Steinhaus's Theorem);

\noindent (ii) $\Sigma $ locally non-meagre at all points $x\in \Sigma $ (by
the Piccard-Pettis Theorem -- such sets can be extracted as subsets of a
second-category set, using separability or by reference to the Banach
Category Theorem);

\noindent (iii) $\Sigma $ universally measurable and not Haar null at any
point (by the Christensen-Solecki Interior-points Theorem -- again such sets
can be extracted using separability).

If $\Sigma $ has the Baire property and is locally non-meagre then it is
co-meagre (since its quasi interior is eveywhere dense).

\bigskip

For contrast with Corollary 2 below, we first note that local boundedness of
mid-point convex functions follows from regularity almost exactly as in the
subadditive case of Theorem 2 of \S 3.

\bigskip

\noindent \textbf{Theorem 2}$^{\prime }$. \textit{For }$X$ \textit{a vector
space with a Steinhaus-Weil, Baire topology }$\mathcal{T}$\textit{\ and} $%
S:X\rightarrow \mathbb{R}$\textit{\ mid-point convex:} \textit{if} $S$
\textit{is }$\mathcal{T}$\textit{-Baire, then it is locally bounded.}

\bigskip

\noindent \textbf{Proof.} Suppose $|S(u+z_{n})|\rightarrow \infty $ for some
$u\in X$ and null sequence $z_{n}\rightarrow 0.$ As the level sets $%
H_{n}^{\pm }:=\{t:|S(\pm t)|\leq n\}$ are $\mathcal{T}$ -Baire and $\mathcal{%
T}$ is Baire, for some $k$ the set $H_{k}^{\pm }$ is non-meagre. As $%
\mathcal{T}$ is Steinhaus-Weil, $H_{k}^{\pm }-H_{k}^{\pm }$ has 0 in its
interior.

First suppose that $S(u+z_{n})\rightarrow +\infty .$ Then there is $n\in
\mathbb{N}$ such that $4z_{m}\in H_{k}^{\pm }-H_{k}^{\pm }$ for all $m\geq n$
$.$ For $m\geq n,$ choose $a_{m},b_{m}\in H_{k}^{\pm }$ with
\[
4z_{m}=a_{m}-b_{m}
\]%
for $m\geq n.$ Then, as%
\[
u+z_{m}=\frac{1}{2}2u+\frac{1}{4}a_{m}+\frac{1}{4}(-b_{m}),
\]%
for all $m\geq n$%
\[
S(u+z_{m})\leq \frac{1}{2}S(2u)+\frac{1}{4}S(a_{m})+\frac{1}{4}S(-b_{m})\leq
\frac{1}{2}S(2u)+\frac{1}{2}k,
\]%
contradicting upper unboundedness.

If on the other hand $S(u+z_{n})\rightarrow -\infty ,$ then argue similarly,
but now choose $k,n$ and $a_{m},b_{m}\in H_{k}^{\pm }$ so that
\[
-2z_{m}=a_{m}-b_{m},
\]%
for all $m\geq n.$ Then%
\[
S(u/2)-\frac{1}{4}S(a_{m})-\frac{1}{4}S(-b_{m})\leq \frac{1}{2}S(u+z_{m}),
\]%
contradicting lower unboundedness. $\square $

\bigskip

This result immediately yields a Banach-space version of Theorem BD in the
separable context. The non-separable variant must wait.

\bigskip

\noindent \textbf{Corollary 1. }\textit{For a separable Banach space }$X$%
\textit{, if} $S:X\rightarrow \mathbb{R}$\textit{\ mid-point convex is Baire
or universally measurable, then it is locally bounded and so continuous.}

\bigskip

\noindent \textbf{Proof. }Apply Theorem 2 or 2H respectively to the $\sigma $%
-ideal of meagre or Haar null sets. $\square $

\bigskip

As with Theorem 2H (at the end of \S 3) so too here, Theorem 2$^{\prime }$
has a `Haar'-type variant with the same proof, which we need below in
Theorem FS.

\bigskip

\noindent \textbf{Theorem 2H}$^{\prime }$. \textit{For }$X$ \textit{a
topological vector space, }$\mathcal{H}$\textit{\ a proper }$\sigma $\textit{%
-ideal with the Steinhaus-Weil property\ and} $S:X\rightarrow \mathbb{R}$%
\textit{\ mid-point convex:} \textit{if} $S$ \textit{is universally
measurable, then it is locally bounded.}

\bigskip

Our aim now is to identify in Theorem 5 below, for any weakly convex
function on $\mathbb{R}$, a canonical continuous convex function using
continuity on sets $\Sigma $ with the local Steinhaus-Weil property.
Thereafter in Theorem 6 we will deduce continuity of a weakly convex
function on $\mathbb{R}$, which we extend to the separable Banach context of
Theorem 7. As corollaries we then deduce Theorems M and FS of \S 1.

\bigskip

\noindent \textbf{Theorem 5 (Canonical Extension Theorem). }\textit{For }$%
\Sigma $\textit{\ locally Steinhaus-Weil and }$I\subseteq \mathrm{%
int(cl(\Sigma ))}$\textit{, if }$S:\mathbb{R}\rightarrow \mathbb{R}$\textit{%
\ is both continuous and mid-point convex on }$\Sigma $, \textit{and}
\[
\bar{S}(x)=\bar{S}^{\Sigma }(x):=\lim \sup\nolimits_{y\rightarrow x}^{\Sigma
}S(y)\qquad (x\in I)
\]%
\textit{-- then}

\noindent (i) \textit{the limit exists for }$x\in I$\textit{: }$\bar{S}%
(x):=\lim_{y\rightarrow x}^{\Sigma }S(y)$\textit{;}

\noindent (ii) $\bar{S}=S$ \textit{on }$\Sigma $\textit{;}

\noindent (iii) \textit{for all }$x\in I$\textit{\ }%
\[
S(x)\leq \bar{S}(x)\text{;}
\]

\noindent (iv)\textit{\ }$\bar{S}$ \textit{is }$\mathbb{R}_{+}$\textit{%
-convex:}%
\[
\bar{S}(tx+(1-t)y)\leq t\bar{S}(x)+(1-t)\bar{S}(y)\qquad (x,y\in I,t\in
(0,1)).
\]

For the proof  we need three lemmas.

\bigskip

\noindent \textbf{Lemma 1 (Full convexity on }$\Sigma $\textbf{). }\textit{%
For} $\Sigma \subseteq \mathbb{R}$\textit{\ locally Steinhaus-Weil, if }$S$
\textit{is both continuous and weakly convex on }$\Sigma $\textit{, then }$S$%
\textit{\ is fully convex on }$\Sigma :$%
\[
S((1-t)a+tb)\leq (1-t)S(a)+tS(b)\qquad (a,b\in \Sigma ,t\in (0,1)).
\]

\bigskip

\noindent \textbf{Proof.} For any $T,$ write $B_{\varepsilon
}^{T}(x):=B_{\varepsilon }(x)\cap T.$ Take any $u.$ Choose $a,b\in \Sigma $
with $a<u<b$ and define $t$ by%
\[
u=(1-t)a+tb:\qquad t=(u-a)/(b-a).
\]%
As $\Sigma -u$ has the Steinhaus-Weil property locally, and exponentiation
is a homeomorphism, for small enough $\varepsilon $%
\[
B_{\varepsilon }^{\Sigma -u}(b-u)[B_{\varepsilon }^{\Sigma
-u}(u-a)]^{-1}+1=-B_{\varepsilon }^{\Sigma -u}(b-u)B_{\varepsilon }^{\Sigma
-u}(a-u)^{-1}+1
\]%
has $(b-u)(u-a)^{-1}+1>1$ in its interior, and so has a rational element $%
r>1.$

Taking successively $\varepsilon =1/n$ for $n\in \mathbb{N},$ select as
above rational $r_{n}>1$ and $a_{n},b_{n}$ in $\Sigma $ such that
\[
a_{n}\rightarrow a,b_{n}\rightarrow b,\text{ }r_{n}=1+\frac{b_{n}-u}{u-a_{n}}%
=\frac{b_{n}-a_{n}}{u-a_{n}}\rightarrow \frac{b-a}{u-a}=1/t.
\]%
So with $q_{n}=1/r_{n}\in \mathbb{Q}_{+},$%
\[
u=a_{n}+q_{n}(b_{n}-a_{n})=(1-q_{n})a_{n}+q_{n}b_{n},\text{ and }0<q_{n}<1.
\]%
As $a,b$ are relative-continuity points and $q_{n}$ is rational with $%
q_{n}\rightarrow t,$%
\begin{eqnarray*}
S(u) &=&S((1-q_{n})a_{n}+q_{n}b_{n}) \\
&\leq &(1-q_{n})S(a_{n})+q_{n}S(b_{n})\rightarrow (1-t)S(a)+tS(b).
\end{eqnarray*}%
So for any $a,b\in \Sigma $ and $0<t<1,$%
\[
S((1-t)a+tb)\leq (1-t)S(a)+tS(b).
\]%
That is: $S$ is $\mathbb{R}_{+}$-convex over $\Sigma .$ $\square $

\bigskip

\noindent \textbf{Corollary 2 (Boundedness on }$\Sigma $\textbf{). }\textit{%
For }$\Sigma \subseteq \mathbb{R}$\textit{\ locally Steinhaus-Weil, if }$S$
\textit{is both continuous and mid-point convex on }$\Sigma $\textit{, then
for each }$x\in $ \textrm{int(cl(}$\Sigma $\textrm{))}\textit{\ and each
sequence }$\{u_{n}\}$ \textit{in} $\Sigma $ \textit{converging to }$x$%
\textit{\ the sequence }$\{S(u_{n})\}$ \textit{is bounded}\textrm{.}

\bigskip

\noindent \textbf{Proof. }For $x\in $ \textrm{int(cl(}$\Sigma $\textrm{)), }%
choose $a,b\in \Sigma $ with $a<x<b;$ then $S$ is bounded above on $(a,b).$
Indeed, applying Lemma 1 to $a,b\in \Sigma ,$%
\[
S((a,b))\leq \max (S(a),S(b)).
\]%
Suppose that $S(u_{n})\rightarrow -\infty $ for some $u_{n}$ in $\Sigma $
with $u_{n}\rightarrow x\in $ \textrm{int(cl(}$\Sigma $\textrm{)). }Take $%
x<v\in \Sigma $ and put $w=(x+v)/2.$ Write $w=t_{n}u_{n}+(1-t_{n})v$ for
some $0<t_{n}<1.$ W.l.o.g. $t_{n}$ is convergent, to $t$ say; then $%
w=tx+(1-t)v=(x+v)/2$ and so $t=1/2.$ But%
\[
S(w)=S(t_{n}u_{n}+(1-t_{n})v)\leq t_{n}S(u_{n})+(1-t_{n})S(v),
\]%
giving in the limit $S(w)\leq -\infty ,$ a contradiction. $\square $

\bigskip

The following result is stated as we need it -- for the line; we raise, and
leave open here, the question of whether it holds in an infinite-dimensional
Banach space. It does, however, hold under a stronger $\mathbb{Q}$-convexity
assumption on the set $\Sigma $ -- see Lemma 2$^{\prime }$ below.

\bigskip

\noindent \textbf{Lemma 2 (Unique limits on }$\mathbb{R}$\textbf{). }\textit{%
For }$\Sigma \subseteq \mathbb{R}$\textit{\ locally Steinhaus-Weil, if }$%
S|\Sigma $\ \textit{is both continuous and mid-point convex, then for any }$%
x\in \mathbb{R}$ \textit{and for any sequences in }$\Sigma $ \textit{with }$%
u_{n}\uparrow x$\textit{\ and }$v_{n}\downarrow x,$%
\[
\lim S(u_{n})=\lim S(v_{n}),
\]%
\textit{when both limits exist.}

\bigskip

\noindent \textbf{Proof. }Put $L:=\lim S(u_{n}),R:=\lim S(v_{n});$ we show
that $L=R.$ If not, suppose first that $L<R.$ For $\varepsilon :=(R-L)/3>0$
there is $m(0)$ so that for $n>m(0),$%
\[
R-\varepsilon \leq S(v_{n}).
\]%
Choose $t_{n}\downarrow 0$ and $m(n)>n$ in order to express the right-sided
sequence $v$ in terms of the left-sided sequence $u$:%
\[
v_{m(n)}=(1-t_{n})u_{m(n)}+t_{n}v_{n}.
\]%
This is possible as $u_{n}\uparrow x$ and $v_{n}\downarrow x.$ As $u_{m(n)}$%
, $v_{m(n)}\in \Sigma ,$ by Lemma 1,%
\[
R-\varepsilon \leq S(v_{m(n)})\leq
(1-t_{n})S(u_{m(n)})+t_{n}S(u_{n})\rightarrow L.
\]%
But $R-\varepsilon \leq L$ gives the contradiction $R-L\leq $ $\varepsilon
\leq (R-L)/3.$

Now suppose that $L>R.$ Taking $\varepsilon =(L-R)/3,$ proceed to a similar
contradiction by exchanging the roles of the $u$ and $v$ sequences: $%
u_{m(n)}=(1-t_{n})v_{m(n)}+t_{n}u_{n}$ with $t_{n}\downarrow 0,$ to obtain%
\[
L-\varepsilon \leq S(u_{m(n)})\leq
(1-t_{n})S(v_{m(n)})+t_{n}S(u_{n})\rightarrow R.
\]%
This time $L-R\leq \varepsilon .$ $\square $

\bigskip

\noindent \textbf{Lemma 2}$^{\prime }$\textbf{\ (Banach-space variant of
unique limits). }\textit{For a Banach space }$X$ \textit{and }$w\in X$%
\textit{, and }$\mathbb{Q}$\textit{-convex (closed under rational convex
combinations) }$\Sigma $\textit{, if }$S:X\rightarrow \mathbb{R}$\textit{\
is both mid-point convex, and locally bounded on }$\Sigma $\textit{\ at }$w$%
\textit{, then for any sequences }$u_{n}\rightarrow x$\textit{\ and }$%
v_{n}\rightarrow x$ \textit{in} $\Sigma $ \textit{with }$\{S(u_{n})\}$
\textit{and }$\{S(v_{n})\}$ \textit{convergent}%
\[
\lim S(u_{n})=\lim S(v_{n}).
\]

\bigskip

\noindent \textbf{Proof. }Put $A:=\lim S(u_{n}),B:=\lim S(v_{n}).$ By
symmetry of the assumptions we may assume that $A<B.$ Noting that the
translate $w+\Sigma $ is $\mathbb{Q}$-convex and the translate $%
S_{w}(x)=S(a+x)$ is mid-point convex on $w+\Sigma ,$ w.l.o.g suppose that $%
x=0.$ Choose $\delta >0$ and $K$ such that $|S(y)|\leq K$ for all $y\in
\Sigma $ with $|y|\leq \delta .$ For $\varepsilon :=(B-A)/3>0,$ there is $%
m(0)$ so that for $n>m(0),$%
\[
B-\varepsilon \leq S(v_{n}).
\]%
Let $t_{n}\downarrow 0$ be dyadic rational, e.g. $t_{n}=2^{-n}.$ Then $%
s_{n}:=1/t_{n}\rightarrow \infty .$ For each $n$ choose $m(n)>n$ such that $%
||u_{m(n)}||<\delta /3$ and $||s_{n}u_{m(n)}||<\delta
/3,||s_{n}v_{m(n)}||<\delta /3.$ Put%
\[
w_{n}:=s_{n}v_{m(n)}+(s_{n}-1)u_{m(n)}\in \Sigma .
\]%
Then%
\[
||w_{n}||=||s_{n}v_{m(n)}||+||s_{n}u_{m(n)}||+||u_{m(n)}||\leq \delta ,
\]%
so that $||S(w_{n})||<K$ and%
\[
v_{m(n)}=(1-t_{n})u_{m(n)}+t_{n}w_{n}.
\]%
Here $S$ is mid-point convex on $\Sigma $, so%
\[
R-\varepsilon \leq S(v_{m(n)})\leq
(1-t_{n})S(u_{m(n)})+t_{n}S(w_{n})\rightarrow A.
\]%
But $B-\varepsilon \leq A$ gives the contradiction $B-A\leq $ $\varepsilon
\leq (B-A)/3.$ $\square $

\bigskip

Below we write $\lim \sup\nolimits_{y\rightarrow x}^{\Sigma }$, $%
\lim_{y\rightarrow x}^{\Sigma }S(y)$\ for the upper limit or limit of $S(y)$
as $y$\ tends to $x$\ through $\Sigma $.

\bigskip

\noindent \textbf{Lemma 3 (Regularization).} \textit{For }$\Sigma $\textit{\
locally Steinhaus-Weil and }$I\subseteq \mathrm{int(cl(\Sigma ))}$, \textit{%
if }$S:\mathbb{R}\rightarrow \mathbb{R}$\textit{\ with }$S|\Sigma $\textit{\
mid-point convex locally bounded}, \textit{\ write}
\[
\bar{S}(x):=\lim \sup\nolimits_{y\rightarrow x}^{\Sigma }S(y)\qquad (x\in
I),
\]%
\textit{Then}

\noindent (i) \textit{the limit exists for }$x\in I$\textit{: }$\bar{S}%
(x):=\lim_{y\rightarrow x}^{\Sigma }S(y)$\textit{;}

\noindent (ii) \textit{the function }$\bar{S}(x)$ \textit{is continuous on }$%
I$\textit{.}

\bigskip

\noindent \textbf{Proof. }(i) By Lemma 2, $\bar{S}$ is well-defined.

\noindent (ii) Suppose that $\bar{S}(x_{n})\rightarrow L\neq \bar{S}(x)$ for
some sequence $x_{n}\rightarrow x$ in $I,$ with $L$ possibly infinite.
Choose $y_{n}\in B_{1/n}(x_{n})\cap \Sigma $ with $|S(y_{n})-\bar{S}%
(x_{n})|<2^{-n}.$ Then $y_{n}\rightarrow x$ and $S(y_{n})\rightarrow L,$
contradicting $S(y_{n})\rightarrow \bar{S}(x).$ $\square $

\bigskip

\noindent \textbf{Proof of Theorem 5. }By Corollary 2, (i) and (ii) follow
as in the proof of Lemma 3, but in (ii) take $y_{n}\in \Sigma $.

By continuity of $S$ on $\Sigma $, $S|\Sigma =\bar{S}|\Sigma .$

\noindent (iii) As in Lemma 1, for any $x\in I$ take $%
x=t_{x}u_{x}+(1-t_{x})v_{x}$ with $u_{x}<x<v_{x},$ $u_{x},v_{x}\in \Sigma ,$
and $t_{x}\in (0,1);$ then%
\[
S(x)\leq t_{x}S(u_{x})+(1-t_{x})S(v_{x}).
\]%
Taking limits as $u_{x}\uparrow x,v_{x}\downarrow x$, and w.l.o.g. assuming $%
t_{x}\rightarrow \tau _{x}$ (by boundedness),
\[
S(x)\leq \tau _{x}\bar{S}(x)+(1-\tau _{x})\bar{S}(x)=\bar{S}(x).
\]

\noindent (iv) Take $x,y,\alpha $ arbitrary, and put $\beta =1-\alpha .$ In $%
\Sigma $ choose $x_{n}\rightarrow x,y_{n}\rightarrow y$ and $%
z_{n}\rightarrow \alpha x+\beta y,$ so that with $\beta _{n}=1-\alpha _{n}$%
\[
z_{n}:=\alpha _{n}x_{n}+\beta _{n}y_{n}:\qquad \alpha
_{n}:=(y_{n}-z_{n})/(y_{n}-x_{n})\rightarrow \lbrack y-\alpha x+\beta
y]/(y-x)=\alpha .
\]%
Then, as $x_{n},y_{n},z_{n}\in \Sigma ,$ from%
\[
S(\alpha _{n}x_{n}+\beta _{n}y_{n})\leq \alpha _{n}S(x_{n})+\beta
_{n}S(y_{n})
\]%
we get%
\[
\bar{S}(\alpha x+\beta y)\leq \alpha \bar{S}(x)+\beta \bar{S}(y).\qquad
\square
\]

\bigskip

For $S:\mathbb{R}\rightarrow \mathbb{R}$ Baire/measurable, since $S$ is
quasi $\sigma $-continuous (Th. BL, \S 2), there is $\Sigma
=\bigcup\nolimits_{m}\Sigma _{m},$ which is quasi all of $\mathbb{R}$ and so
dense in $\mathbb{R}$, with $\Sigma _{m}$ an increasing sequence of sets
\textit{each having the Steinhaus-Weil property locally}. Before returning
to a Banach space setting we prove a result in $\mathbb{R}$, which we shall
apply (twice) later in the context of a `typical' line segment.

\bigskip

\noindent \textbf{Theorem 6. }\textit{For a dense set }$\Sigma
=\bigcup\nolimits_{m}\Sigma _{m}$\textit{, with each }$\Sigma _{m}$\ \textit{%
locally Steinhaus-Weil},\textit{\ if }$S:\mathbb{R}\rightarrow \mathbb{R}$%
\textit{\ is mid-point convex on }$\Sigma $ \textit{and quasi-}$\sigma $-%
\textit{continuous with respect to }$\Sigma $\textit{, then }$S\ $\textit{is
continuous.}

\bigskip

\noindent \textbf{Proof. }If not, referring to the dense set $\Sigma
=\bigcup\nolimits_{m}\Sigma _{m}$ and the continuous functions $\bar{S}%
^{\Sigma _{m}}$ of the Canonical Extension Theorem, which, identified with
their graphs, form an increasing union (by (i) above), we may put $\bar{S}%
:=\bigcup\nolimits_{m}\bar{S}^{\Sigma _{m}}$ and so may suppose w.l.o.g. for
some $x>0$ that $S(x)<\bar{S}(x).$ (Otherwise shift the origin to the left,
or consider $S(-x).)$ Fix such an $x,$ and put $\varepsilon :=[\bar{S}%
(x)-S(x)]/4.$ By continuity of $\bar{S}$ at $x,$ for some $\Delta \in (0,x),$%
\begin{equation}
|\bar{S}(x)-\bar{S}(y)|\leq \varepsilon \qquad (y\in (x-\Delta ,x+\Delta )).
\tag{*}
\end{equation}%
Take $0<\delta <(\Delta /2),$ and set $z:=x+\delta .$

For some $m$, $\Sigma ^{\prime }:=\Sigma _{m}\cap (x,z)$ is non-empty and so
has the Steinhaus-Weil property; then $\Sigma ^{\prime }+\Sigma ^{\prime }$
contains an interval, $(a,b)$ say. Put $s:=\inf \Sigma ^{\prime }\geq x;$
then $(a,b)\subseteq (2s,2z)\subseteq (2x,2z).$ As $a\geq 2s,$ there is $%
t\in \Sigma ^{\prime }$ with $t>s,$ such that $\alpha t\in (a,b)$ for some
dyadic rational $\alpha >2.$ Indeed, for $2s<a$ and $t\in (s,a/2)\cap \Sigma
^{\prime },$ taking $\alpha \in (a/t,b/t)$ gives $2t<a<\alpha t<b$ and $%
\alpha >2;$ on the other hand, for $2s=a,$ if $t\in (s,b/2)\cap \Sigma
^{\prime },$ then $2<b/t,$ and so taking $\alpha \in (2,b/t),$ gives $%
2s=a<2t<\alpha t<b.$

For any dyadic $\alpha >2,$ take%
\[
q=q(\alpha )=\alpha -1>1;
\]%
then $q$ is a positive dyadic rational, and%
\[
\frac{1}{\alpha }+\frac{q}{\alpha }=1.
\]%
Furthermore, if $\alpha $ satisfies $\alpha t<2z,$ then as $t>x$%
\[
1<q=\alpha -1<\frac{z+z-x}{x}=\frac{x+2(z-x)}{x}=1+2\frac{\delta }{x}.
\]%
For $t>s$ in $\Sigma ^{\prime }$ and dyadic $\alpha >2$ as above, since $%
\alpha t\in \Sigma ^{\prime }+\Sigma ^{\prime },$
\[
\alpha t=(x+u)+(x+v),\text{ }\qquad t=(1/\alpha )x+(q/\alpha )[x+u+v]/q,
\]%
with $0<u,z<\delta .$ As $\delta <(\Delta /2)/(1-\Delta /x)$ (as $1-\Delta
/x<1),$ as above.%
\[
1<q<1+2\frac{\delta }{x}<\frac{\Delta }{x-\Delta }+1<x/(x-\Delta ).
\]%
So%
\[
(x-\Delta )<x/q<x.
\]%
Also as $\delta <\Delta /2$ and $1<q,$%
\[
(u+v)/q<2\delta /q<2\delta <\Delta .
\]%
So for $\delta <\Delta /2$ small enough as above, $[x+u+v]/q\in \lbrack
x-\Delta ,x+\Delta ]$ and likewise $(x+u)/q\in \lbrack x-\Delta ,x+\Delta ].$
As $S(y)\leq \bar{S}(y),$ for all $y,$ and as $t\in \Sigma _{m},$ using
mid-point convexity and continuity of $\bar{S}$ at $t$ (indeed of $\bar{S}%
^{\Sigma _{m}}$)$,$%
\begin{eqnarray*}
\bar{S}(x)-\varepsilon &\leq &\bar{S}(t)=S(t)\leq (1/\alpha )S(x)+(q/\alpha
)S([x+u+v]/q) \\
&\leq &(1/\alpha )S(x)+(q/\alpha )\bar{S}([x+u+v]/q) \\
&\leq &(1/\alpha )S(x)+(q/\alpha )[\bar{S}(x)+\varepsilon ].
\end{eqnarray*}%
So for $\delta >0$ small enough%
\[
\bar{S}(x)-\varepsilon \leq (1/\alpha )S(x)+(q/\alpha )[\bar{S}%
(x)+\varepsilon ],
\]%
where $\alpha $ and $q$ depend on $\delta .$ But
\[
1\leq \frac{1}{\alpha }(1+1+2\delta /x)\text{,}\qquad \frac{1}{2}>\frac{1}{%
\alpha }\geq \frac{1}{2(1+\delta /x)}.
\]%
Let $\delta \downarrow 0:$ $1/\alpha \rightarrow 1/2,$ so
\[
\bar{S}(x)-\varepsilon \leq (1/2)S(x)+(1/2)[\bar{S}(x)+\varepsilon ],
\]%
or%
\[
\bar{S}(x)-S(x)\leq 3\varepsilon =(3/4)[\bar{S}(x)-S(x)],
\]%
a contradiction. $\square $

\bigskip

As a corollary we now have a result on separable Banach spaces, which by
Theorem B will enable us to prove in their more general setting Theorems M
and FS, stated in \S 1. Note the local character of the key assumption.

\bigskip

\noindent \textbf{Theorem 7. }\textit{For a separable Banach space }$X,$
\textit{a dense set }$\Sigma =\bigcup\nolimits_{m}\Sigma _{m}$\textit{, with
each }$\Sigma _{m}$\ \textit{locally Steinhaus-Weil},\textit{\ if }$%
S:X\rightarrow \mathbb{R}$\textit{\ is mid-point convex on }$\Sigma ,$
\textit{Baire, and quasi-}$\sigma $-\textit{continuous with respect to }$%
\Sigma $\textit{,} \textit{then }$S\ $\textit{is continuous.}

\bigskip

\noindent \textbf{Proof. }Since $\Sigma $ has the Steinhaus-Weil property
locally, we may proceed as in Theorem 6 above to consider $x\neq 0$ with $%
S(x)<\bar{S}(x);$ define $\varepsilon >0$ as there and choose $\Delta >0$
similarly so that (*) holds for $y\in B_{\Delta }(x).$ Take $\delta <\Delta
/2$ and $\Sigma ^{\prime }:=\Sigma \cap B_{\delta }(x).$ By the
Kuratowski-Ulam Theorem, for some $\sigma \in \Sigma ^{\prime }$ the ray
\[
R_{x}(\sigma ):=\{x+\lambda (\sigma -x):\lambda \geq 0\}
\]%
meets $\Sigma ^{\prime }$ in a non-meagre set: otherwise $\Sigma ^{\prime
}\cap R_{x}(\sigma )$ is meagre for all $\sigma \in \Sigma ^{\prime },$ and
so $\Sigma ^{\prime }$ is meagre. As $\Sigma ^{\prime }\cap R_{x}(\sigma )$
is Baire there is an interval $I:=[s,s^{\prime }]$ along $R_{x}(\sigma )$
for which $\Sigma ^{\prime }\cap I$ is co-meagre in $I.$ Continue as in
Theorem 6 working in $R_{x}(\sigma )$ rather than $\mathbb{R}_{+}$ to obtain
a contradiction to $S(x)<\bar{S}(x),$ so deducing continuity of $S$. $%
\square $

\bigskip

As an immediate corollary we are now able to prove Theorem M due to Mehdi
(albeit for a general topological vector space), and Theorem FS, a result
slightly weaker than of Fischer and S\l odkowski [FisS] (where universal
measurability is modulo Haar null sets).

\bigskip

\noindent \textbf{Proof of Theorem M. }By Theorem B\ we may assume w.l.o.g.
that $X$ is separable. By Theorem BL, $S$ is continuous relative to a
co-meagre (so dense) set $\Sigma .$ Since $\Sigma $ has the Steinhaus-Weil
property locally, we may apply Theorem 7 above with $\Sigma _{m}\equiv
\Sigma $, as $S$ is mid-point convex on $\Sigma ,$ so deducing continuity of
$S$. $\square $

\bigskip

\noindent \textbf{Proof of Theorem FS. }As above, we may again assume that $%
X $ is separable. For any distinct points $a,b,$ consider the line $L$
through $a$ and $b,$ and let $\lambda $ be Lebesgue masure on $L.$ Then $%
S|L:L\rightarrow \mathbb{R}$ is universally measurable, so $\lambda $%
-measurable and so quasi-$\sigma $-continuous by Luzin's Theorem. By Theorem
6, $S|L$ is continuous on $L$ and so fully convex on $L$. So $S\ $is fully
convex. By Theorem 2H$^{\prime }$, $S$ is locally bounded, so continuous by
Theorem BD*. $\square $

\bigskip

We close with an analogue of Theorem 7. We will need to argue as in Theorem
6 twice: once, in the `measure-case' mode of Theorem 6 (using $\sigma $%
-continuity), to establish that the continuity points form a big set (as in
Luzin's Theorem), and then again, but now in the `category mode' of Theorem
6 as in Theorem 7 (where $\Sigma $ is dense and locally Steinhaus-Weil).
This reflects the hybrid nature of Christensen's definition of Haar null
sets.

\bigskip

\noindent \textbf{Theorem 8. }\textit{For a separable Banach space }$X,$
\textit{a dense set }$\Sigma =\bigcup\nolimits_{m}\Sigma _{m}$\textit{, with
each }$\Sigma _{m}$\ \textit{locally Steinhaus-Weil},\textit{\ if }$%
S:X\rightarrow \mathbb{R}$\textit{\ is mid-point convex on }$\Sigma $
\textit{and universally measurable,} \textit{then }$S\ $\textit{is
continuous.}

\bigskip

\noindent \textbf{Proof. }Put $\Gamma :=\{x\in X:S$ is continuous at $x\}$;
then $\Gamma $ is universally measurable. Indeed, by Lemma 3 $\bar{S}$ is
well-defined and continuous (from the given $\Sigma $). Thus $S$ is
discontinuous at $x$ iff $S(x)\neq \bar{S}(x),$ and so, since $\bar{S}$ is
continuous and $S$ universally measurable, the complement of $\Gamma $ is%
\begin{eqnarray*}
\bigcup\nolimits_{q\in \mathbb{Q}}\{x &:&S(x)<q<\bar{S}(x)\}\cup \{x:\bar{S}%
(x)<q<S(x)\} \\
&=&\bigcup\nolimits_{q\in \mathbb{Q}}S^{-1}(-\infty ,q)\cap \bar{S}%
^{-1}(q,\infty )\cup \bar{S}^{-1}(-\infty ,q)\cap S^{-1}(q,\infty ),
\end{eqnarray*}%
so universally measurable.

We claim first that $\Gamma \cap U$ is non-Haar null for all non-empty open $%
U.$ If not, $U\cap \Gamma $ is Haar null for some non-empty open $U$; then,
by the definition of Haar nullity (see \S 3), there exist a Borel set $%
G\supseteq U\cap \Gamma $ and a Borel probability measure $\mu $ such that $%
\mu (g+G)=0$ for all $g\in X.$ W.l.o.g. $U=u+B_{\delta };$ as $X$ is
separable, a countable number of translates $t_{i}+U$ of $U,$ and so also of
$B_{\delta },$ covers $X$. So $\mu (u+v+B_{\delta })>0$ for some $v:=t_{i}.$
Put $\mu _{\nu }(E)=\mu (v+E)$ for $E\subseteq X$ Borel; then $\mu _{v}$ is
finite with $\mu _{v}(U)>0$, and $S$ is quasi-$\sigma $-continuous w.r.t. $%
\mu _{v},$ by Luzin's Theorem. Proceed as in Theorem 7, but this time
applying Christensen's WFT in place of the Kuratowski-Ulam Theorem (again
since $S$ is universally measurable), to deduce that $S$ is continuous at $x$
for each $x\in U$, so contradicting the assumption that $G$ is Haar null
(and so not the whole of $B_{\delta }$).

Being universally measurable and locally non-Haar null, $\Gamma $ has the
Steinhaus-Weil property locally, by a theorem of Christensen [Chr1, Th. 2]
(extended by Solecki [Sole3, Th. 1(ii) via Prop. 3.3(i)]). With $\Sigma
=\Gamma $ and $X=\bar{\Gamma}$, proceed once more as in Theorem 7, again
applying Christensen's WFT in place of the Kuratowski-Ulam Theorem. This
gives that $S$ is continuous on $X$. $\square $

\section{Complements}

\noindent \textit{1. Berz's other theorems. }A sublinear function $S$ has $%
\mathbb{Q}_{+}$-convex epigraph $C$. This observation allows Berz to deduce
from the $\mathbb{Q}$-version of the Hahn-Banach theorem that $S$ is the
supremum of all the additive functions $f$ which it majorizes; the proof
refers to the $\mathbb{Q}$-hyperplanes defined by $f$ that support the
epigraph. Since a Baire/measurable $S$ is locally bounded (Th. 2 above), all
of the additive minorants of $S$ supporting $C$ are bounded above and so
linear by \textit{Darboux's Theorem} (see e.g. [BinO9] and the references
cited there). This allows Berz to deduce that their upper envelope comprises
the two half-lines defining $S$ (equivalently, this is the upper envelope of
the supremum and infimum of the additive minorants of $S$). Hence Berz
deduces a third result: when $S$ is symmetric about the origin it may be
represented as a \textit{norm}. Indeed, embed $x\mapsto \{f(x)\}_{f}$ so
that $f(x)$ is the projection of $x$ onto the $f$ co-ordinate space; then a
norm is defined by
\[
||x||:=\sup_{f}|f(x)|=S(x).
\]

\noindent \textit{2. Automatic continuity. }The proof of Theorem 1 is
inspired by an idea due to Goldie appearing in [BinG, I, Th. 5.7] (cf.
[BinGT, Th. 3.2.5]), and more fully exploited in a recent series of papers
including [BinO10-12, 13 Prop. 3]. The theme here is the interplay between
\textit{functional inequalities} (as with subadditivity, convexity etc.) and
\textit{functional equations} (as with additivity and the Cauchy functional
equation). Here, minimal regularity implies continuity -- whence the term
automatic continuity -- and linearity; see e.g. [BinO8] and the references
cited there.

\noindent \textit{3. Automatic continuity and group action. }An automatic
continuity theorem of Hoffmann-J\o rgensen is particularly relevant here for
the discussion of the Baire-Berz Theorem. Hoffmann-J\o rgensen proves in
[Rog, Part 3: Th. 2.2.12] the (sequential) continuity of a Baire function $%
f:X\rightarrow Y$ when a single non-meagre group $T$ acts on the two
(Hausdorff) spaces $X$ and $Y$ with $f(tx)=tf(x),$ by appealing to a KBD
argument (under $T$ rather than under addition) in $X.$ In the Baire-Berz
Theorem it is a meagre group, namely $\mathbb{Q}_{+}$, that acts
multiplicatively on the Banach spaces $X$ and $Y=\mathbb{R}$; but it is the
additive structure of a Banach space which permits the use of KBD to obtain
global continuity from continuity on a smaller set.

\noindent \textit{4. Convex and coherent risk measures.} As remarked in \S %
1, Berz's sublinearity theorem is connected with the theory of coherent risk
measures [FolS, \S\ 4.1]. The key properties are convexity and positive
homogeneity ($\rho (\lambda x)=\lambda \rho (x)$ for $\lambda \geq 0$).
Under positive homogeneity, convexity is equivalent to subadditivity. This
paper thus extends to sublinearity studies of the related areas of
convexity, subadditivity and additivity, for which see e.g. [BinO1, 3].

In the economic/financial context, positive homogeneity -- a form of
scale-invariance -- means that large and small firms (or agents) have
similar preferences; see e.g. Lindley [Lin, Ch. 5]. This is far from the
case in practice, which is why convex risk measures (in which positive
homogeneity is dropped) are often preferred; again, see e.g. [FolS, \S 4.1].
Sensitivity to scale here is related to curvature of utility functions, and
the `law of diminishing returns'. This incidentally underpins the viability
of the insurance industry; again see e.g. Lindley [Lin, Ch. 5].

The two half-lines in Berz's theorem correspond to taking long and short
positions in one dimension. One can extend to many dimensions, as in [FolS],
where the `broken line' becomes a cone, and as we do in \S 4. Berz himself
worked in one dimension, as his motivation was normability (below).

\noindent \textit{5. Normability. }As norms are necessarily sublinear,
Berz's third result (6.1) addresses the question of which sublinear
functions are realized as norms. In this connection, the criterion for
normability of a topological vector space was established by Kolmogorov, see
e.g. [Rud, Th. 1.39]; for recent metric characterizations of normability --
in terms of translation-invariant metrics -- see the Oikhberg-Rosenthal
result [OikR] demanding continuity of scaling and \textit{isometry} of all
one-dimensional subspaces $R(x)$ with $\mathbb{R}$. \v{S}emrl's relaxation
[Sem] drops this continuity when spaces are of dimension at least $2$. (As
for relaxation of homogeneity see [Mat].) Invariant metrics are provided by
the Birkhoff-Kakutani normability theorem -- see e.g. [Rud, Th. 1.24],
[HewR, Th. 8.3], or for recent accounts [Gao, Ch. 1-4], [Ost2, \S 2.1].

\noindent 6. \textit{Beyond local compactness:\ Haar category-measure
duality. }In the absence of Haar measure, the definition (in \S 2) of left
Haar null subsets of a topological group $G$ required $\mathcal{U}(G),$ the
universally measurable sets -- by dint of the role of the totality of
(probability) measures on $G$. The natural dual of $\mathcal{U}(G)$ is the
class $\mathcal{U}_{\mathcal{B}}(G)$ of \textit{universally Baire sets},
defined,for $G$ with a Baire topology, as those sets $B$ whose preimages $%
f^{-1}(B)$ are Baire (have the Baire property) in any compact Hausdorff
space $K$ for any continuous $f:K\rightarrow G$. Initially considered in
[FenMW] for $G=\mathbb{R}$, these have attracted continued attention for
their role in the investigation of axioms of determinacy and large cardinals
-- see especially [Woo]; cf. [MarS].

Analogously to the left Haar null sets, define in $G$ the family of \textit{%
left Haar meagre} sets, $\mathcal{HM}(G)$, to comprise the sets $M$
coverable by a universally Baire set $B$ for which there are a compact
Hausdorff space $K$ and a continuous $f:K\rightarrow G$ with $f^{-1}(gB)$
meagre in $K$ for all $g\in G.$ These were introduced, in the \textit{abelian%
} Polish group setting and with $K$ metrizable, by Darji [Dar], cf. [Jab],
and shown there to form a $\sigma $-ideal of meagre sets (co-extensive with
the meagre sets for $G\ $locally compact); as $\mathcal{HM}(G)\mathcal{%
\subseteq B}_{0}(G),\mathcal{\ }$the family is not studied here.

\newpage

\bigskip

\noindent \textbf{References.}

\noindent \lbrack Ajt] M. Ajtai, On the boundedness of definable linear
operators, \textsl{Period. Math. Hungar.} \textbf{5} (1974), 343--352.%
\newline
\noindent \lbrack ArtDEH] P. Artzner, F. Delbaen, J.-M. Eber, D. Heath,
Coherent measures of risk. \textsl{Math. Finance} \textbf{9} (1999), no. 3,
203--228.\newline
\noindent \lbrack BalR] M. Balcerzak, J. Rzepecka, Marczewski sets in the
Hashimoto topologies for measure and category. \textsl{Acta Univ. Carolin.
Math. Phys. }\textbf{39 }(1998), 93-97\newline
\noindent \lbrack BarFN] A. Bartoszewicz, Artur, M. Filipczak, T. Natkaniec,
On Smital properties. \textsl{Topology Appl.} \textbf{158} (2011),
2066--2075.\newline
\noindent \lbrack Ber] E. Berz, Sublinear functions on $\mathbb{R},$ \textrm{%
\textsl{Aequat. Math.}} \textbf{12} (1975), 200-206.\newline
\noindent \lbrack BinG] N. H. Bingham, C.M. Goldie, Extensions of regular
variation. I. Uniformity and quantifiers; Extensions of regular variation.
II. Representations and indices. \textsl{Proc. London Math. Soc.} (3)
\textbf{44} (1982), 473--496, 497--534.\newline
\noindent \lbrack BinGT] N. H. Bingham, C. M. Goldie and J. L. Teugels,
\textsl{Regular variation}, 2nd ed., Cambridge University Press, 1989 (1st
ed. 1987). \newline
\noindent \lbrack BinO1] N. H. Bingham and A. J. Ostaszewski, Generic
subadditive functions, \textsl{Proc. Amer. Math. Soc. }\textbf{136} (2008),
4257-4266.\newline
\noindent \lbrack BinO2] N. H. Bingham and A. J. Ostaszewski, Infinite
combinatorics and the foundations of regular variation, \textsl{Journal of
Math. Anal. Appl.} \textbf{360} (2009), 518-529.\newline
\noindent \lbrack BinO3] N. H. Bingham and A. J. Ostaszewski, Automatic
continuity: subadditivity, convexity, uniformity. \textsl{Aequationes
Mathematicae}, 78 (2009), 257-270\newline
\noindent \lbrack BinO4] N. H. Bingham and A. J. Ostaszewski, Beyond
Lebesgue and Baire: generic regular variation, \textsl{Coll. Math.} \textbf{%
116} (2009), 119-138.\newline
\noindent \lbrack BinO5] N. H. Bingham and A. J. Ostaszewski, Kingman,
category and combinatorics. \textsl{Probability and Mathematical Genetics}
(Sir John Kingman Festschrift, ed. N. H. Bingham and C. M. Goldie), 135-168,
London Math. Soc. Lecture Notes in Mathematics \textbf{378}, CUP, 2010.
\newline
\noindent \lbrack BinO6] N. H. Bingham and A. J. Ostaszewski, Beyond
Lebesgue and Baire II: Bitopology and measure-category duality. \textsl{%
Colloquium Math.} \textbf{121} (2010), 225-238.\newline
\noindent \lbrack BinO7] N. H. Bingham and A. J. Ostaszewski, Regular
variation without limits. \textsl{J. Math. Anal. Appl.}, \textbf{370}
(2010), 322-338.\newline
\noindent \lbrack BinO8] N. H. Bingham and A. J. Ostaszewski, Normed versus
topological groups: dichotomy and duality, \textsl{Dissertationes Math.}
\textbf{472} (2010) 138 pp.\newline
\noindent \lbrack BinO9] N. H. Bingham and A. J. Ostaszewski, Dichotomy and
infinite combinatorics: the theorems of Steinhaus and Ostrowski. \textsl{%
Math. Proc. Camb. Phil. Soc.} \textbf{150} (2011), 1-22. \newline
\noindent \lbrack BinO10] N. H. Bingham and A. J. Ostaszewski, The Steinhaus
theorem and regular variation : De Bruijn and after, \textsl{Indag. Math.}
\textbf{24} (2013),679-692. \newline
\noindent \lbrack BinO11] N. H. Bingham and A. J. Ostaszewski, Cauchy's
functional equation and extensions: Goldie's equation and inequality, the Go%
\l \k{a}b-Schinzel equation and Beurling's equation, \textsl{Aequat. Math.}
\textbf{89} (2015), 1293-1310 (fuller version: arXiv: 1405.3947).\newline
\noindent \lbrack BinO12] N. H. Bingham and A. J. Ostaszewski, Beurling
moving averages and approximate homomorphisms, \textsl{Indag. Math. (NS) }
\textbf{27} (2016), 601--633 (fuller version: arXiv1407.4093v2).\newline
\noindent \lbrack BinO13] N. H. Bingham and A. J. Ostaszewski, Additivity,
subadditivity and linearity: automatic continuity and quantifier weakening,
arXiv: 1405.3948v2.\newline
\noindent \lbrack BinO14] N. H. Bingham and A. J. Ostaszewski, Beyond
Lebesgue and Baire IV: Density topologies and a converse Steinhaus-Weil
Theorem, arXiv1607.00031.\newline
\noindent \lbrack BinO15] N. H. Bingham and A. J. Ostaszewski, The
Steinhaus-Weil property: its converse, Solecki amenability and
subcontinuity, arXiv1607.00049.\newline
\noindent \lbrack Blu] H. Blumberg, On convex functions, \textsl{Trans.
Amer. Math. Soc.} \textbf{20} (1919), 40-44.\newline
\noindent \lbrack BorD] D. Borwein and S. Z. Ditor, Translates of sequences
in sets of positive measure, \textsl{Canadian Math. Bull.} \textbf{21}
(1978), 497-498.\newline
\noindent \lbrack BorM] J. M. Borwein, W. B. Moors, Null sets and
essentially smooth Lipschitz functions. \textsl{SIAM J. Optim.}\textbf{\ 8}%
.2 (1998), 309--323.\newline
\noindent \lbrack Bru] A. M. Bruckner, Differentiation of integrals. \textsl{%
Amer. Math. Monthly} \textbf{78}.9 (1971), 1-51\newline
\noindent \lbrack Chr1] J. P. R. Christensen, On sets of Haar measure zero
in abelian Polish groups. Proceedings of the International Symposium on
Partial Differential Equations and the Geometry of Normed Linear Spaces
(Jerusalem, 1972). \textsl{Israel J. Math.} \textbf{13} (1972), 255--260
(1973).\newline
\noindent \lbrack Chr2] J. P. R. Christensen, \textsl{Topology and Borel
structure. Descriptive topology and set theory with applications to
functional analysis and measure theory.} North-Holland Mathematics Studies,
Vol. 10, 1974.\newline
\noindent \lbrack FenN] J. E. Fenstad, D. Normann, On absolutely measurable
sets. \textsl{Fund. Math.} \textbf{81} (1973/74), no. 2, 91--98.\newline
\noindent \lbrack FisS] P. Fischer, Z. S\l odkowski, Christensen zero sets
and measurable convex functions.\textsl{\ Proc. Amer. Math. Soc.} \textbf{79}%
.3 (1980), 449--453.\newline
\noindent \lbrack FolS] H. F\"{o}llmer and A. Schied, \textsl{Stochastic
finance: An introduction in discrete time}, de Gruyter Studies in Math.
\textbf{27}, 3$^{\text{rd}}$ ed., 2011 (1$^{\text{ed}}$ 2002, 2$^{\text{nd}}$
ed. 2004).\newline
\noindent \lbrack Fre1] D. Fremlin, Measure-additive coverings and
measurable selectors, \textsl{Dissertationes Math.} \textbf{260} (1987),116p.%
\newline
\noindent \lbrack Fre2] D. Fremlin, \textsl{Measure theory Vol. 4:
Topological measure theory.} Torres-Fremlin, 2$^{\text{nd}}$ printing, 2006
(1$^{\text{st}}$ printing 2003).\newline
\noindent \lbrack FreNR] D. Fremlin, T. Natkaniec, I. Rec\l aw, Universally
Kuratowski-Ulam spaces, \textsl{Fund. Math.} 165 (2000), no. 3, 239--247.%
\newline
\noindent \lbrack Gao] S. Gao, \textsl{Invariant descriptive theory}, CRC\
Press, 2009.\newline
\noindent \lbrack Gar] H. G. Garnir, Solovay's axiom and functional
analysis, 189--204, in \textsl{Functional analysis and its applications},
eds. H. G. Garni, K. R. Unni, J.H. Willaimson, Proc. Madras Conf. on Funct.
Anal., Lecture Notes in Math. \textbf{399}, Springer, 1974.\newline
\noindent \lbrack GofNN] C. Goffman, C. J. Neugebauer and T. Nishiura,
Density topology and approximate continuity, \textsl{Duke Math. J.} \textbf{%
28} (1961), 497--505.\newline
\noindent \lbrack GofW] C. Goffman, D. Waterman, Approximately continuous
transformations, \textsl{Proc. Amer. Math. Soc.} \textbf{12} (1961),
116--121.\newline
\noindent \lbrack GroE1] K.-G. Grosse-Erdmann, Regularity properties of
functional equations and inequalities, \textsl{Aequat. Math.} \textbf{37}
(1989) 233-251.\newline
\noindent \lbrack GroE2] K.-G. Grosse-Erdmann, An extension of the
Steinhaus-Weil Theorem, \textsl{Coll. Math.} \textbf{57} (1989), 307-317.%
\newline
\noindent \lbrack Hal] P. R. Halmos, \textsl{Measure Theory}, Grad. Texts in
Math. \textbf{18}, Springer 1974 (1st. ed. Van Nostrand, 1950).\newline
\noindent \lbrack HarLP] G. H. Hardy, J. E. Littlewood, G. P\'{o}lya,
\textsl{Inequalities}. 2nd ed., CUP, 1952.\newline
\noindent \lbrack Has] H. Hashimoto, On the *topology and its application.
\textsl{Fund. Math. }\textbf{91} (1976), no. 1, 5--10.\newline
\noindent \lbrack HewR] E. Hewitt, K. A. Ross, \textsl{Abstract harmonic
analysis,} \textsl{Vol. I. Structure of topological groups, integration
theory, group representations}. 2$^{\text{nd}}$ ed. Grundlehren Wiss., 115.
Springer, 1979.\newline
\noindent \lbrack HilP] E. Hille and R. S. Phillips, \textsl{Functional
analysis and semi-groups}, Coll. Publ. Vol \textbf{31}, Amer. Math. Soc, 3rd
ed. 1974 (1st ed. 1957).\newline
\noindent \lbrack HunSY] B.R. Hunt, T. Sauer, J. A. Yorke, Prevalnce: a
translation-invariant \textquotedblleft almost everywhere\textquotedblright\
on infinite dimensional spaces, \textsl{Bull. Amer. Math. Soc.} \textbf{27}%
.2 (1992), 217-238.\newline
\noindent \lbrack Jab] E. Jab\l o\'{n}ska, Some analogies between Haar
meager sets and Haar null sets in abelian Polish groups. \textsl{J. Math.
Anal. Appl.} \textbf{421} (2015),1479--1486.\newline
\noindent \lbrack JanH] D. Jankovi\'{c}, T. R. Hamlett, New topologies from
old via ideals, \textsl{Amer. Math. Monthly}, \textbf{97} (1990), 295-310.%
\newline
\noindent \lbrack JayR] J. E. Jayne and C. A. Rogers, \textsl{Selectors.}
Princeton University Press, 2002.\newline
\noindent \lbrack KalPR] N. J. Kalton, N. T. Peck, J. W. Roberts, \textsl{An
F-space sampler.} London Math. Soc. Lect. Note Ser. \textbf{89}, Cambridge
University Press, 1984.\newline
\noindent \lbrack Kec] A. S. Kechris: \textsl{Classical Descriptive Set
Theory.} Grad. Texts in Math. \textbf{156}, Springer, 1995.\newline
\noindent \lbrack Kem] J. H. B. Kemperman, A general functional equation,
\textsl{Trans. Amer. Math. Soc.} \textbf{86} (1957), 28--56.\newline
\noindent \lbrack Kes] H. Kestelman, The convergent sequences belonging to a
set, \textsl{J. London Math. Soc.\ }\textbf{22} (1947), 130-136.\newline
\noindent \lbrack Kom] Z. Kominek, On an equivalent form of a Steinhaus
theorem, \textsl{Mathematica (Cluj)}, \textbf{30 (53)} (1988), 25-27.\newline
\noindent \lbrack Kuc] M. Kuczma, \textsl{An introduction to the theory of
functional equations and inequalities. Cauchy's equation and Jensen's
inequality.} 2nd ed., Birkh\"{a}user, 2009 [1st ed. PWN, Warszawa, 1985].%
\newline
\noindent \lbrack Lin] D. V. Lindley, \textsl{Taking decisions}, 2nd ed.,
Wiley, 1985.\newline
\noindent \lbrack LukMZ] J. Luke\v{s}, J. Mal\'{y}, L. Zaj\'{\i}\v{c}ek,
\textsl{Fine topology methods in real analysis and potential theory},
Lecture Notes in Mathematics \textbf{1189}, Springer, 1986.\newline
\noindent \lbrack Mat] J. Matkowski, Subadditive functions and a relaxation
of the homogeneity condition of seminorms. \textsl{Proc. Amer. Math. Soc. }%
\textbf{117} (1993), no. 4, 991--1001\newline
\noindent \lbrack McNFE] A. J. McNeil, R. Frey, P. Embrechts, \textsl{%
Quantitative risk management. Concepts, techniques and tools}. Princeton
University Press, 2005.\newline
\noindent \lbrack Meh] M. R. Mehdi, On convex functions. \textsl{J. London
Math. Soc. }\textbf{39} (1964), 321--326.\newline
\noindent \lbrack vMil] J. van Mill, A note on the Effros Theorem, \textsl{%
Amer. Math. Monthly}, 111.9 (2004), 801-806.\newline
\noindent \lbrack vMilP] J. van Mill, R. Pol, The Baire category theorem in
products of linear spaces and topological groups$\QTR{sl}{.}$ \textsl{%
Topology Appl.}. 22 (1986), no. 3, 267--282.\newline
\noindent \lbrack Mil] H. I. Miller, Generalization of a result of Borwein
and Ditor, \textsl{Proc. Amer. Math. Soc.} \textbf{105} (1989), no. 4,
889--893.\newline
\noindent \lbrack MilO] H. I. Miller and A. J. Ostaszewski, Group action and
shift-compactness, \textsl{J. Math. Anal. App.} \textbf{392} (2012), 23-39.%
\newline
\noindent \lbrack OikR] T. Oikhberg, H. Rosenthal, A metric characterization
of normed linear spaces, \textsl{Rocky Mountain J. Math.} \textbf{37} (2007)
597--608.\newline
\noindent \lbrack Ost1] A. J. Ostaszewski, Analytically heavy spaces:
Analytic Cantor and Analytic Baire Theorems, \textsl{Topology and its
Applications} \textbf{158} (2011), 253-275.\newline
\noindent \lbrack Ost2] A. J. Ostaszewski, Beyond Lebesgue and Baire III:
Steinhaus's Theorem and its descendants, \textsl{Topology and its
Applications} \textbf{160} (2013), 1144-1154.\newline
\noindent \lbrack Ost3] A. J. Ostaszewski, Almost completeness and the
Effros Theorem in normed groups, \textsl{Topology Proceedings} \textbf{41}
(2013), 99-110 (fuller version: arXiv.1606.04496).\newline
\noindent \lbrack Ost4] A. J. Ostaszewski, Shift-compactness in almost
analytic submetrizable Baire groups and spaces, survey article, \textsl{%
Topology Proceedings} \textbf{41} (2013), 123-151.\newline
\noindent \lbrack Ost5] A. J. Ostaszewski, Effros, Baire, Steinhaus and
non-separability, \textsl{Topology and its Applications}, Mary Ellen Rudin
Special Issue, \textbf{195 }(2015), 265-274.\newline
\noindent \lbrack Oxt] J. C. Oxtoby: \textsl{Measure and category}, 2nd ed.
Graduate Texts in Math. \textbf{2}, Springer, 1980 (1$^{\text{st}}$ ed.
1972).\newline
\noindent \lbrack Par] K. R. Parthasarathy, \textsl{Probability measures on
metric spaces.} Academic Press, 1967.\newline
\noindent \lbrack Pol] R. Pol, Note on category in Cartesian products of
metrizable spaces, \textsl{Fund. Math.} 102 (1979), no. 1, 55--59.\newline
\noindent \lbrack Roc] R. T. Rockafellar, \textsl{Convex analysis},
Princeton, 1972.\newline
\noindent \lbrack Rog] C. A. Rogers, J. Jayne, C. Dellacherie, F. Tops\o e,
J. Hoffmann-J\o rgensen, D. A. Martin, A. S. Kechris, A. H. Stone, \textsl{%
Analytic sets.} Academic Press, 1980.\newline
\noindent \lbrack Ros] R. A. Rosenbaum, Sub-additive functions, \textsl{Duke
Math. J.} \textbf{17} (1950), 227-242.\newline
\noindent \lbrack RosS] C. Rosendal, S. Solecki, Automatic continuity of
homomorphisms and fixed points on metric compacta. \textsl{Israel J. Math.}
\textbf{162} (2007), 349--371\newline
\noindent \lbrack Rud] W. Rudin, \textsl{Functional Analysis}, McGraw-Hill,
1991.\newline
\noindent \lbrack Sem] P. \v{S}emrl, A characterization of normed spaces.
\textsl{J. Math. Anal. Appl.} \textbf{343} (2008), 1047--1051.\newline
\noindent \lbrack Sie1] W. Sierpi\'{n}ski, Sur l'\'{e}quation fonctionelle%
\textsl{\ }$f(x+y)=f(x)+f(y)$, \textsl{Fund. Math.} \textbf{1} (1920),
116-122 (reprinted in \textsl{Oeuvres choisis,} Vol. 2, 331-33, PWN, Warsaw
1975).\newline
\noindent \lbrack Sie2] W. Sierpi\'{n}ski, Sur les fonctions convexes
mesurables, \textsl{Fund. Math.} \textbf{1} (1920), 125-128 (reprinted in
\textsl{Oeuvres choisis,} Vol. 2, 337-340, PWN, Warsaw 1975).\newline
\noindent \lbrack Sim] B. Simon, \textsl{Convexity: An analytic viewpoint},
Cambridge Tracts in Math. \textbf{187}, CUP, 2011.\newline
\noindent \lbrack Sole1] S. Solecki, On Haar null sets. \textsl{Fund. Math. }%
\textbf{149} (1996), no. 3, 205--210.\newline
\noindent \lbrack Sole2] S. Solecki, Size of subsets of groups and Haar null
sets, \textsl{Geom. Funct. Analysis}, \textbf{15} (2005), 246-273.\newline
\noindent \lbrack Sole3] S. Solecki, Amenability, free subgroups, and Haar
null sets in non-locally compact groups. \textsl{Proc. London Math. Soc.}
\textbf{(3) 93} (2006), 693--722.\newline
\noindent \lbrack Sole4] S. Solecki, A Fubini Theorem. \textsl{Topology
Appl. }\textbf{154} (2007), no. 12, 2462--2464.\newline
\noindent \lbrack Trau] R. Trautner, A covering principle in real analysis,
\textsl{Quart. J. Math. Oxford} (2) \textbf{38} (1987), 127- 130.\newline
\noindent \lbrack Wei] A. Weil, \textsl{L'integration dans les groupes
topologiques}, Actualit\'{e}s Scientifiques et Industrielles 1145, Hermann,
1965 (1$^{\text{st }}$ ed. 1940).\newline
\noindent \lbrack Wil] W. Wilczy\'{n}ski, \textsl{Density topologies},
Handbook of measure theory, Vol. I, II, 675--702, North-Holland, 2002.%
\newline
\noindent \lbrack Wri1] J. D. Maitland Wright, On the continuity of
mid-point convex functions. \textsl{Bull. London Math. Soc.}\textbf{\ 7}
(1975), 89--92.\newline
\noindent \lbrack Zak] E. Zakon, A remark on the theorems of Lusin and
Egoroff. \textsl{Canad. Math. Bull.} \textbf{7} (1964), 291--295.\newline

\bigskip

\bigskip

\noindent Mathematics Department, Imperial College, London SW7 2AZ;
n.bingham@ic.ac.uk \newline
Mathematics Department, London School of Economics, Houghton Street, London
WC2A 2AE; A.J.Ostaszewski@lse.ac.uk\newpage

\section*{Appendix 1: Set-theoretic foundations\protect\footnote{%
This Appendix for the arXiv version only.}}

We summarize below background information needed to appreciate the various
set-theoretic axioms to which we have referred. As mentioned in the
Introduction, this may be omitted by the expert (or uninterested) reader;
the earlier article [Wri2] of 1977 had a similar motivation.

\bigskip

\noindent \textit{1. Category/measure regularity versus practicality. }The
Baire/measurable property assumed above is usually satisfied in mathematical
practice. Indeed, any analytic subset of $\mathbb{R}$ possesses these
properties ([Rog, Part 1 \S 2.9], [Kec, 29.5]), hence so do all the sets in
the $\sigma $-algebra that they generate (the $C$-sets, [Kec, \S 29.D]).
There is a broader class still. Recall first that an analytic set may be
viewed as a projection of a planar Borel set $P,$ so is definable as $%
\{x:\Phi (x)\}$ via the $\mathbf{\Sigma }_{1}^{1}$ formula $\Phi
(x):=(\exists y\in \mathbb{R})[(x,y)\in P];$ here the notation $\mathbf{%
\Sigma }_{1}^{1}$ indicates one quantifier block (the subscripted value) of
existential quantification, ranging over reals (type 1 objects -- the
superscripted value). Use of the \textit{bold-face version} of the symbol
indicates the need to refer to \textit{arbitrary} coding (by reals not
necessarily in an \textit{effective} manner, for which see [Gao, \S 1.5])
the various opens sets needed to construct $P.$ (An open set $U$ is \textit{%
coded} by the sequence of rational intervals contained in $U.)$

Consider a set $A$ such that both $A$ and $\mathbb{R}\backslash A$ may be
defined by a $\mathbf{\Sigma }_{2}^{1}$ formula, say respectively as $%
\{x:\Phi (x)\}$ and $\{x:\Psi (x)\}$, where $\Phi (x):=(\exists y\in \mathbb{%
R})(\forall z\in \mathbb{R})(x,y,z)\in P\}$ now, and similarly $\Psi .$ This
means that $A$ is both $\mathbf{\Sigma }_{1}^{1}$ and $\mathbf{\Pi }_{1}^{1}$
(with $\mathbf{\Pi }$ indicating a leading universal quantifier block), and
so in the ambiguous class $\mathbf{\Delta }_{1}^{1}.$ If in addition the
equivalence%
\[
\Phi (x)\Longleftrightarrow \lnot \Psi (x)
\]%
is provable in $ZF,$ i.e. \textit{without reference} to $AC,$ then $A$ is
said to be \textit{provably} $\mathbf{\Delta }_{2}^{1}.$ It turns out that
such sets have the Baire/measurable property -- see [FenN], where these are
generalized to the \textit{universally (=absolutely) measurable }sets (\S %
2). How much further this may go depends on what axioms of set-theory are
admitted, a matter to which we now turn.

Our interest in such matters is dictated by the \textit{Character Theorems}
of regular variation, noted in [BinO7, \S 3] (revisited in [BinO12]), which
identify the logical complexity of the function $h^{\ast }(x):=\lim \sup
h(t+x)-h(t),$ which is $\mathbf{\Delta }_{2}^{1}$ if the function $h$ is
Borel (and is $\mathbf{\Pi }_{2}^{1}$ if $h$ is analytic, and $\mathbf{\Pi }%
_{3}^{1}$ if $h$ is co-analytic). We argued in [BinO7, \S 5] that $\mathbf{%
\Delta }_{2}^{1}$ is a natural setting in which to study regular variation.

\noindent \textit{2. Principle of D\textit{ependent Choice DC}. }In his
paper Berz relied on the \textit{Axiom of Choice} AC, in the usual form of
\textit{Zorn's Lemma},\ which is used in the same context of $\mathbb{R}$
over the field of scalars $\mathbb{Q}$ as in Hamel's construction of a
discontinuous additive function, and so ultimately rests on transfinite
induction of \textit{continuum} length requiring \textit{continuum} many
selections. Our proof of Berz's theorem depends in effect on the Baire
Category Theorem BC, or the completeness of $\mathbb{R}$, since Theorem KBD
is a variant of BC (see \S 2), and so ultimately rests on elementary
induction via the \textit{Axiom (Principle) of Dependent Choice(s)} DC (thus
named in 1948 by Tarski [Tar2, p. 96] and studied in [Mos], but anticipated
in 1942 by Bernays' [Ber, Axiom IV*, p. 86] -- see [Jec1, \S 8.1], [Jec2,
Ch. 5]), and DC is equivalent to BC by a result of Blair [Bla]. (For further
results in this direction see also [Pin1,2], [Gol], [HerK], [Wol], and the
textbook [Her].)

We note that DC is equivalent to a statement about trees: a pruned tree has
an infinite branch (for which see [Kec, 20.B]) and so by its very nature is
an ingredient in set-theory axiom systems which consider the extent to which
Banach-Mazur-type games (with underlying tree structure) are determined. The
latter in turn have been viewed as generalizations of Baire's Theorem ever
since Choquet [Cho] -- cf. [Kec, 8C,D,E]. Inevitably, determinacy and the
study of the relationship between category and measure go hand in hand.

\noindent \textit{3. Practical axiomatic alternatives: LM, PB, AD, PD. }%
While ZF is common ground in mathematics, AC is not, and alternatives to it
are widely used, in which for example all sets are Lebesgue-measurable
(usually abbreviated to LM) and all sets have the Baire property, sometimes
abbreviated to PB (as distinct from BP to indicate individual `possession of
the Baire property'). One such DC above. As Solovay [Solo2, p. 25] points
out, this axiom is sufficient for the establishment of Lebesgue measure,
i.e. including its translation invariance and countable additivity
("...positive results ... of measure theory..."), and may be assumed
together with LM. Another is the \textit{Axiom of Determinacy} AD mentioned
above and introduced by Mycielski and Steinhaus [MycS]; this implies LM, for
which see [MySw], and PB, the latter a result due to Banach -- see [Kec,
38.B]. Its introduction inspired remarkable and still current developments
in set theory concerned with determinacy of `definable' sets of reals (see
[ForK] and particularly [Nee]) and consequent combinatorial properties (such
as partition relations) of the alephs (see [Kle]). Others include the
(weaker) \textit{Axiom of Projective Determinacy} PD [Kec, \S\ 38.B],
restricting the operation of AD\ to the smaller class of projective sets.
(The independence and consistency of DC\ versus AD was established
respectively in Solovay [Solo3] and Kechris [Kech] -- see also [KechS].)

\noindent \textit{4. LM versus PB. }In 1983 Raissonier and Stern [RaiS, Th.
2] (cf. [Bar1,2]), inspired by then current work of Shelah (circulating in
manuscript since 1980) and earlier work of Solovay, showed that if every $%
\mathbf{\Sigma }_{2}^{1}$ set is Lebesgue measurable, then every $\mathbf{%
\Sigma }_{2}^{1}$ set has BP, whereas the converse fails -- for the latter
see [Ste] -- cf. [BarJ, \S 9.3]. This demonstrates that measurability is in
fact the stronger notion -- see [JudSh, \S 1] for a discussion of the
consistency of analogues at level 3 and beyond -- which is why we regard
category rather than measure as \textit{primary}. For we have seen above how
the category version of Berz's theorem implies its measure version; see also
[BinO6,7].

Note that the assumption of G\"{o}del's \textit{Axiom of Constructibility} $%
V=L,$ a strengthening of AC, yields $\mathbf{\Delta }_{2}^{1}$
non-measurable subsets, so that the Fenstad-Normann result on the narrower
class of provably $\mathbf{\Delta }_{2}^{1}$ sets mentioned in 6.1 above
marks the limit of such results in a purely ZF framework (at level 2).

\noindent \textit{5. Consistency and the role of large cardinals. }While LM\
and PB are inconsistent with AC, such axioms can be consistent with DC.
Justification with scant exception involves some form of large-cardinal
assumption, which in turn calibrates relative consistency strengths -- see
[Kan] and [KoeW] (cf. [Lar] and [KanM]). Thus Solovay [Solo2] in 1970 was
the first to show the equiconsistency of ZF+DC+LM+PB with that of ZFC+`%
\textit{there exists an inaccessible cardinal}'. The appearance of the
inaccessible in this result is not altogether incongruous, given its
emergence in results (from 1930 onwards) due to Banach [Ban] (under GCH),
Ulam [Ula] (under AC), and Tarski [Tar1], concerning the cardinalities of
sets supporting a countably additive/finitely additive [0,1]-valued/$\{0,1\}$%
-valued measure (cf. [Bog, 1.12(x)]. Later in 1984 Shelah [She1, 5.1] showed
in ZF+DC that already the measurability of all $\mathbf{\Sigma }_{1}^{3}$
sets implies that $\aleph _{1}^{L}$ is inaccessible (the symbol $\aleph
_{1}^{L}$ refers to the substructure of constructible sets and denotes the
first uncountable ordinal therein). As a consequence, Shelah [She1,5.1A]
showed that ZF+DC+LM\ is equiconsistent with ZF+`\textit{there exists an
inaccessible}', whereas [She1, 7.17] ZF+DC+PB\ is equiconsistent with just
ZFC (i.e. without reference to inaccessible cardinals), so driving another
wedge between classical measure-category symmetries (see [JudSh] for
further, related `wedges'). The latter consistency theorem relies on the
result [She1, 7.16] that any model of ZFC\ + CH has a generic (forcing)
extension satisfying ZF+ `\textit{every set of reals (first-order) defined
using a real and an ordinal parameter has BP}'. For a topological proof see
Stern [Ste].

\noindent \textit{6. LM versus PB continued. }Raisonnier [Rai, Th. 5] (cf.
[She1, 5.1B]) has shown that in ZF+DC one can prove that if there is an
uncountable well-ordered set of reals (in particular a set of cardinality $%
\aleph _{1}$), then there is a non-measurable set of reals. (This motivates
Judah and Spinas [JudSp] to consider generalizations including the
consistency of the $\omega _{1}$-variant of DC.) See also Judah and Ros\l %
anowski [JudR] for a model (due to Shelah) in which ZF+DC+LM+$\lnot $PB\
holds, and also [She2] where an inaccessible cardinal is used to show
consistency of ZF+LM+$\lnot $PB+`\textit{there is an uncountable set without
a perfect subset}'. For a textbook treatment of much of this material see
again [BarJ].

Raisonnier [Rai, Th. 3] notes the result, due to Shelah and Stern, that
there is a model for ZF+DC+PB+$\aleph _{1}=\aleph _{1}^{L}$+ `\textit{the
ordinally definable subsets of real are measurable}'. So, in particular by
Raisonnier's result, there is a non-measurable set in this model. Shelah's
result indicates that the non-measurable is either $\Sigma _{3}^{1}$
(light-face symbol: all open sets coded effectively) or $\mathbf{\Sigma }%
_{2}^{1}$ (bold-face). Thus here PB+$\lnot $LM holds.

\noindent \textit{7. Regularity of reasonably definable sets. }From the
existence of suitably large cardinals flows a most remarkable result due to
Shelah and Woodin [SheW] justifying the opening practical remark about BP,
which is that every `reasonably definable' set of reals is Lebesgue
measurable: compare the commentary in [BecK] following their Th 5.3.2. This
is a latter-day sweeping generalization of a theorem due to Solovay (cf.
[Solo1]) that, subject to large-cardinal assumptions, $\mathbf{\Sigma }%
_{2}^{1}$ sets are measurable (and so also have BP by [RaiS]).

\bigskip

\noindent \textbf{References for Appendix 1.}

\noindent \lbrack Bar1] T. Bartoszy\'{n}ski, Additivity of measure implies
additivity of category. \textsl{Trans. Amer. Math. Soc.} \textbf{281}
(1984), no. 1, 209--213 \newline
\noindent \lbrack Bar2] T. Bartoszy\'{n}ski, \textsl{Invariants of measure
and category}, Ch. 7 in [ForK].\newline
\noindent \lbrack BarJ] T. Bartoszy\'{n}ski and H. Judah, \textsl{Set
Theory: On the structure of the real line}, Peters 1995.\newline
\noindent \lbrack BecK] H. Becker, A. S. Kechris, \textsl{The descriptive
set theory of Polish group actions.} London Math. Soc. Lect. Note Ser.,
\textbf{232}, Cambridge University Press, 1996.\newline
\noindent \lbrack Ber] P. Bernays, A system of axiomatic set theory. III.
Infinity and enumerability. Analysis, \textsl{J. Symbolic Logic} 7 (1942),
65--89.\newline
\noindent \lbrack Bla] C. E. Blair, The Baire category theorem implies the
principle of dependent choices. \textsl{Bull. Acad. Polon. Sci. S\'{e}r.
Sci. Math. Astronom. Phys.} \textbf{25} (1977), 933--934.\newline
\noindent \lbrack Cho] G. Choquet, \textsl{Lectures on analysis}, Vol. I,
Benjamin, New York, 1969.\newline
\noindent \lbrack ForK] M. Foreman, A. Kanamori, \textsl{Handbook of Set
theory}, Springer, 2010.\newline
\noindent \lbrack Fre3] D. Fremlin, \textsl{Measure theory Vol. 5:
Set-theoretic measure theory, Part I.} Torres-Fremlin, 2008.\newline
\noindent \lbrack Gol] R. Goldblatt, On the role of the Baire category
theorem and dependent choice in the foundations of logic. \textsl{J.
Symbolic Logic} \textbf{50} (1985), no. 2, 412--422. \newline
\noindent \lbrack Her] H. Herrlich, \textsl{Axiom of choice.} Lecture Notes
in Mathematics, 1876. Springer, 2006.\newline
\noindent \lbrack HerK] H. Herrlich, K. Keremedis, The Baire category
theorem, and the axiom of dependent choice.\textsl{\ Comment. Math. Univ.
Carolin.} \textbf{40} (1999), no. 4, 771---775.\newline
\noindent \lbrack HowR] P. Howard, J. E. Rubin, The Boolean prime ideal
theorem plus countable choice do [does] not imply dependent choice. \textsl{%
Math. Logic Quart.} \textbf{42} (1996), no. 3, 410--420.\newline
\noindent \lbrack Jac] S. Jackson, \textsl{Structural consequences of AD},
Ch. 21 in [ForK].\newline
\noindent \lbrack Jec1] T. J. Jech, \textsl{The axiom of choice.} Studies in
Logic and the Foundations of Mathematics, Vol. \textbf{75}. North-Holland,
1973.\newline
\noindent \lbrack Jec2] T. J. Jech, \textsl{Set Theory}, 3$^{\text{rd}}$
Millennium ed. Springer, 2003.\newline
\noindent \lbrack JudR] H. Judah, A. Ros\l anowski, On Shelah's
amalgamation. \textsl{Set theory of the reals} (Ramat Gan, 1991), 385--414,
Israel Math. Conf. Proc., \textbf{6}, Bar-Ilan Univ., Ramat Gan, 1993.%
\newline
\noindent \lbrack JudSh] H. Judah, S. Shelah, Baire property and axiom of
choice. \textsl{Israel J. Math.} \textbf{84} (1993), no. 3, 435--450.\newline
\noindent \lbrack JudSp] H. Judah, O. Spinas, Large cardinals and projective
sets. \textsl{Arch. Math. Logic} \textbf{36} (1997), no. 2, 137--155.
\newline
\noindent \lbrack Kan] A. Kanamori, \textsl{The higher infinity. Large
cardinals in set theory from their beginnings}, Springer, 2$^{\text{nd}}$
ed. 2003 (1$^{\text{st}}$ ed. 1994).\newline
\noindent \lbrack KanM] A. Kanamori, M. Magidor, The evolution of large
cardinal axioms in set theory. \textsl{Higher set theory} (Proc. Conf.,
Math. Forschungsinst., Oberwolfach, 1977), pp. 99--275, Lecture Notes in
Math. \textbf{669}, Springer, 1978.\newline
\noindent \lbrack Kech] A. S. Kechris,The axiom of determinacy implies
dependent choices in $L(\mathbb{R})$. \textsl{J. Symbolic Logic} \textbf{49}
(1984), no. 1, 161--173.\newline
\noindent \lbrack Kle] E. M. Kleinberg, \textsl{Infinitary combinatorics and
the axiom of determinateness.} Lecture Notes in Math. \textbf{612},
Springer, 1977.\newline
\noindent \lbrack KechS] A. S. Kechris, R. M. Solovay, On the relative
consistency strength of determinacy hypotheses. \textsl{Trans. Amer. Math.
Soc.} \textbf{290} (1985), no. 1, 179--211.\newline
\noindent \lbrack KoeW] P. Koellner, W.H. Woodin, \textsl{Large cardinals
from determinacy}, Ch. 23 in [ForK].\newline
\noindent \lbrack Lar] P. B. Larson, A brief history of determinacy, \textsl{%
The Cabal Seminar Vol. 4} (eds. A. S. Kechris, B. L\"{o}we, J. R. Steel),
Assoc. Symbolic Logic, 2010.

\noindent \lbrack Mos] Y. Moschovakis, \textsl{Notes on set theory.} 2$^{%
\text{nd}}$ ed. Undergrad. Texts in Math. Springer, 2006.\newline
\noindent \lbrack Most] A. Mostowski, On the principle of dependent choices.
\textsl{Fund. Math.} \textbf{35} (1948). 127--130.\newline
\noindent \lbrack MycS] J. Mycielski, H. Steinhaus, A mathematical axiom
contradicting the axiom of choice. \textsl{Bull. Acad. Polon. Sci. S\'{e}r.
Sci. Math. Astronom. Phys.} \textbf{10} (1962), 1--3.\newline
\noindent \lbrack MySw] J. Mycielski, S. \'{S}wierczkowski, On the Lebesgue
measurability and the axiom of determinateness. \textsl{Fund. Math.} \textbf{%
54} (1964), 67--71.\newline
\noindent \lbrack Nee] I. Neeman, \textsl{Determinacy in }$L(\mathbb{R})$,
Ch. 21 in [ForK].\newline
\noindent \lbrack Pin1] D. Pincus, Adding dependent choice to the prime
ideal theorem. Logic Colloquium 76 (Oxford, 1976), pp. 547--565. \textsl{%
Studies in Logic and Found. Math.}, Vol. \textbf{87}, North-Holland, 1977.%
\newline
\noindent \lbrack Pin2] D. Pincus, Adding dependent choice. \textsl{Ann.
Math. Logic} \textbf{11} (1977), no. 1, 105--145.\newline
\noindent \lbrack Rai] J. Raisonnier, A mathematical proof of S. Shelah's
theorem on the measure problem and related results. \textsl{Israel J. Math.}
\textbf{48} (1984), no. 1, 48--56.\newline
\noindent \lbrack RaiS] J. Raisonnier, J. Stern, Mesurabilit\'{e} and propri%
\'{e}t\'{e} de Baire, \textsl{Comptes Rendus Acad. Sci. I. (Math.)} \textbf{%
296} (1983), 323-326.\newline
\noindent \lbrack She1] S. Shelah, Can you take Solovay's inaccessible away?
\textsl{Israel J. Math.} \textbf{48} (1984), 1-47.\newline
\noindent \lbrack She2] S. Shelah, On measure and category. \textsl{Israel
J. Math.} \textbf{52} (1985), no. 1-2, 110--114.\newline
\noindent \lbrack SheW] S. Shelah, H. Woodin, Large cardinals imply that
every reasonably definable set of reals is Lebesgue measurable. \textsl{%
Israel J. Math.} \textbf{70} (1990), no. 3, 381--394.\newline
\noindent \lbrack Solo1] R. M. Solovay, On the cardinality of $\Sigma
_{2}^{1}$ sets of reals. 1969 Foundations of Mathematics (Symposium
Commemorating Kurt G\"{o}del, Columbus, Ohio, 1966) pp. 58--73 Springer, New
York\newline
\noindent \lbrack Solo2] R. M. Solovay, A model of set-theory in which every
set of reals is Lebesgue measurable. \textsl{Ann. of Math.} (2) \textbf{92}
(1970), 1--56.\newline
\noindent \lbrack Solo3] R. M. Solovay, The independence of DC from AD. The
Cabal Seminar 76--77 (Proc. Caltech-UCLA Logic Sem., 1976--77), pp.
171--183, \textsl{Lecture Notes in Math.} \textbf{689}, Springer, 1978.%
\newline
\noindent \lbrack Ste] J. Stern, Regularity properties of definable sets of
reals, \textsl{Annals of Pure and Appl. Logic}, \textbf{29} (1985), 289-324.%
\newline
\noindent \lbrack Tar1] A. Tarski, Une contribution \`{a} la th\'{e}orie de
la mesure, \textsl{Fund. Math.} \textbf{15} (1930), 42-50.\newline
\noindent \lbrack Tar2] A. Tarski, Axiomatic and algebraic aspects of two
theorems on sums of cardinals. \textsl{Fund. Math.} \textbf{35} (1948),
79--104.\newline
\noindent \lbrack Ula] S. Ulam, Zur Masstheorie in der allgemeinen
Mengenlehre. \textsl{Fund. Math.} \textbf{16} (1930), 140-150.\newline
\noindent \lbrack Wol] E. Wolk, On the principle of dependent choices and
some forms of Zorn's lemma, \textsl{Canad. Bull. Math.}, \textbf{26} (1983),
365-367.\newline
\noindent \lbrack Wri2] J. D. Maitland Wright, Functional Analysis for the
practical man, 283--290 in Functional Analysis: Surveys and Recent Results
\textbf{27}, North-Holland Math. Studies, 1977.\newline
\newpage

\section*{Appendix 2: Blumberg Dichotomy\protect\footnote{%
This Appendix for the arXiv version only.}}

\noindent \textbf{Proof of Theorem B. }The first assertion is established
[Blu, Th. 1] as the second step of an argument but without depending on the
(local compactness) assumptions of the first step (note also that the
quantity $h$ there could be $+\infty $). We expand Blumberg's rather brief
proof below.

As for the second assertion, if $S$ is not continuous at $x_{0}$ choose a
sequence $x_{n}$ with limit $x_{0}$ and with $S(x_{n})$ unbounded, and take $%
D$ to be the closed span of $\{x_{n}:n=0,1,2,...\}.$ Continuity of $S|D$ at $%
x_{0}$ implies that $S(x_{n})$ is bounded, a contradiction.

\bigskip

To return to the first assertion: fix $x_{0}$ with $S$ not continuous at $%
x_{0}.$ First construct a sequence $x_{n}$ convergent to $x_{0}$ with%
\begin{equation}
S(x_{0})<\lim S(x_{n})\leq \infty ,  \tag{$\ddag $}
\end{equation}%
as follows. Begin with any sequence $u_{n}$ with limit $x_{0}$ such that $%
S(u_{n})$ fails to converge to $S(x_{0}).$ The construction now splits into
two cases.

\noindent Case (i): $\lim \sup S(u_{n})>S(x_{0}).$ Passage to a subsequence $%
x_{n}$ of $u_{n}$ yields the desired result that $\lim S(x_{n})>S(x_{0})$.

\noindent Case (ii) $\lim \sup S(u_{n})\leq S(x_{0}).$ Then $\lim \inf
S(u_{n})<S(x_{0}).$ Passing to a subsequence, we may assume that $\lim
S(u_{n})<S(x_{0});$ then taking%
\[
y_{n}:=2x_{0}-u_{n}:\qquad x_{0}=(u_{n}+y_{n})/2
\]%
gives
\[
2S(x_{0})-S(u_{n})\leq S(y_{n}),
\]%
implying
\[
S(x_{0})<2S(x_{0})-\lim S(u_{n})\leq \lim \inf S(y_{n}).
\]%
Now pass to a subsequence $x_{n}$ of $y_{n}$ to obtain $S(x_{0})<\lim
S(x_{n}).$

In either case we obtain a sequence $x_{n}$ with limit $x_{0}$ and with ($%
\ddag $).

Put $h_{n}:=S(x_{n})-S(x_{0});$ then $h:=\lim h_{n}\in (0,\infty ].$
Consider the positive ray from $x_{0}$ to $x_{n}$%
\[
R_{+}(x_{n}):=\{\lambda (x_{0}-x_{n}):\lambda >0\};
\]%
on this ray, for any $k\in \mathbb{N},$ choose $k+1$ equally spaced points,
denoted $x_{n}(i)$ for $i=0,1,...,k,$ starting with $x_{n}(0)=x_{0}$ and $%
x_{n}(1)=x_{n}$ (so that the distance apart of consecutive points is $%
||x_{n}-x_{0}||).$ As above, since
\[
x_{n}(i+1)=\frac{1}{2}(x_{n}(i)+x_{n}(i+2)),\qquad i=0,1,...,k-2
\]%
\[
2S(x_{n}(i+1))\leq S(x_{n}(i))+S(x_{n}(i+2)):\qquad
S(x_{n}(i+1))-S(x_{n}(i))\leq S(x_{n}(i+2))-S(x_{n}(i+1)),
\]%
and so inductively:%
\[
h_{n}=S(x_{n})-S(x_{0})=S(x_{n}(1))-S(x_{n}(0))\leq ...\leq
S(x_{n}(k))-S(x_{n}(k-1)).
\]%
So, using telescoping sums,
\begin{eqnarray*}
S(x_{n}(k))-S(x_{0})
&=&[S(x_{n}(k))-S(x_{n}(k-1))]+[S(x_{n}(k-1))-S(x_{n}(k-2))]+ \\
&&...+S(x_{n})-S(x_{0}) \\
&\geq &kh_{n}.
\end{eqnarray*}%
Taking successively $k=k_{m}:=m,$ and choosing $n=n(m)$ so large that $%
k_{m}||x_{n(m)}-x_{0}||<1/m,$ we obtain a subsequence $x_{n(m)}$ with%
\[
||x_{n(m)}-x_{0}||=k_{m}||x_{n(m)}-x_{0}||\rightarrow 0,
\]%
and, since $h_{n(m)}\rightarrow h\in (0,\infty ],$%
\[
S(x_{n}(k))=S(x_{n(m)}(k_{m}))\geq S(x_{0})+k_{m}h_{m(n)}\rightarrow \infty
.\qquad \square
\]

\end{document}